\documentclass{amsart}
\usepackage{amssymb, latexsym, amsmath}

\usepackage[all]{xy}

\usepackage{enumerate}

\usepackage[usenames, dvipsnames]{color}

\usepackage{mathrsfs} 

\usepackage[OT2,T1]{fontenc}
\DeclareSymbolFont{cyrletters}{OT2}{wncyr}{m}{n}
\DeclareMathSymbol{\Sha}{\mathalpha}{cyrletters}{"58}

\theoremstyle{plain}
\newtheorem{Theorem}[equation]{Theorem}
\newtheorem{Lemma}[equation]{Lemma}
\newtheorem{Corollary}[equation]{Corollary}

\newtheorem*{Conjecture*}{Conjecture}
\newtheorem{Proposition}[equation]{Proposition}

\theoremstyle{definition}
\newtheorem{Definition}[equation]{Definition}

\newtheorem{Paragraph}[equation]{}

\theoremstyle{remark}
\newtheorem{Remark}[equation]{Remark}

\numberwithin{equation}{section}

\renewcommand{\b}{{\text{\rm b}}}
\newcommand{\h}{{\text{\rm h}}}
\renewcommand{\k}{{\kappa}}

\newcommand{\D}{{\mathscr D}}
\newcommand{\E}{{\mathscr E}}

\newcommand{\G}{{\text{\rm G}}}

\newcommand{\Gm}{{\text{\rm G}_{\text{\rm m}}}}
\renewcommand{\H}{{\text{\rm H}}}

\newcommand{\nL}{{{}_n\text{\rm L}}}
\newcommand{\M}{{\text{\rm M}}}
\newcommand{\N}{{\mathbb N}}
\renewcommand{\O}{{\text{\rm O}}}

\renewcommand{\P}{{\mathbb P}}
\newcommand{\Q}{{\mathbb Q}}

\renewcommand{\S}{{\mathcal S}}
\newcommand{\T}{{\mathcal T}}

\newcommand{\cZ}{{\text{\rm Z}}}
\newcommand{\Z}{{\mathbb Z}}

\newcommand{\br}[1]{{\left<{#1}\right>}}

\newcommand{\car}{{\text{\rm char}}}

\newcommand{\cor}{{\text{\rm cor}}}

\newcommand{\cs}{{\text{\rm cs}}}

\newcommand{\df}{{\,\overset{\text{\rm df}}{=}\,}}

\renewcommand{\div}{{\text{\rm div}}}
\renewcommand{\dim}{{\text{\rm dim}}}

\newcommand{\id}{{\text{\rm id}}}
\newcommand{\im}{{\text{\rm im}}}
\newcommand{\ind}{{\text{\rm ind}}}
\renewcommand{\inf}{{\text{\rm inf}\,}}

\newcommand{\isim}{{\;\overset{\sim}{\longrightarrow}\;}}
\newcommand{\isom}{{\;\simeq\;}}
\renewcommand{\ker}{{\text{\rm ker}}}

\newcommand{\nr}{{\text{\rm nr}}}
\newcommand{\ns}{{\text{\rm ns}}}

\newcommand{\ramloc}{{\text{\rm ram.loc.}}}
\newcommand{\red}{{\text{\rm red}}}

\newcommand{\res}{{\text{\rm res}}}
\newcommand{\rk}{{\text{\rm rk}}}

\newcommand{\Ass}{{\text{\rm Ass}}}

\newcommand{\Br}{{\text{\rm Br}}}

\newcommand{\nBrdim}{{{}_n\text{\rm Br.dim}}}

\newcommand{\CaCl}{{\text{\rm CaCl}\,}}

\newcommand{\Div}{{\text{\rm Div}\,}}

\newcommand{\Frac}{{\text{\rm Frac}\,}}

\newcommand{\Length}{{\text{\rm L}}}

\newcommand{\Pic}{{\text{\rm Pic}\,}}

\newcommand{\Supp}{{\text{\rm Supp}}}
\newcommand{\Spec}{{\text{\rm Spec}\,}}

\begin{document}

\title[Cyclic Length in the Brauer Group]
{Cyclic Length in the Tame Brauer Group of the Function Field of a $p$-Adic Curve}

\author{Eric Brussel, Kelly McKinnie, and Eduardo Tengan}
\address 
{Department of Mathematics \& Computer Science\\
Emory University\\
Atlanta, GA 30322\\ USA}
\email{brussel@mathcs.emory.edu}
\address
{Instituto de Ci\^encias Matem\'aticas e de Computa\c c\~ao\\
Universidade de S\~ao Paulo\\
S\~ao Carlos, S\~ao Paulo\\ Brazil}
\email
{etengan@icmc.usp.br}
\address
{University of Montana, Department of Mathematical Sciences, Missoula, MT
  59801, USA}
\email{kelly.mckinnie@umontana.edu}




\begin{abstract}
Let $F$ be the function field of a smooth curve over the $p$-adic number field $\Q_p$.
We show that for each prime-to-$p$ number $n$ the $n$-torsion
subgroup $\H^2(F,\mu_n)={}_n\Br(F)$ is generated by $\Z/n$-cyclic classes; 
in fact the $\Z/n$-length is equal to two.
It follows that the Brauer dimension of $F$ is two (first proved in \cite{Sa97}),
and any $F$-division algebra of period $n$ and index $n^2$ is decomposable.
\end{abstract}

\hfill July 9, 2013

\maketitle

\section*{Introduction}

Let $F$ be a field, $n$ a number that is prime to $\car(F)$,
and let $\H^2(F,\mu_n)$ denote the degree-two Galois cohomology group.
It is well known that $\H^2(F,\mu_n)$ equals ${}_n\Br(F)$,
the $n$-torsion of the Brauer group of $F$, which is by definition the set of equivalence
classes of central simple $F$-algebras of period $n$.
A central simple $F$-algebra is said to be {\it $\Z/n$-cyclic} if it contains a $\Z/n$-cyclic
Galois maximal \'etale subalgebra.
The algebra then has a simple description as a crossed product, and the corresponding
($\Z/n$-cyclic) Brauer class can be written as a cup product $(f)\cdot\theta$ for $(f)\in\H^1(F,\mu_n)$
and $\theta\in\H^1(F,\Z/n)$.
Among all central simple $F$-algebras (resp. Brauer classes), the cyclic $F$-algebras
(resp. cyclic Brauer classes) are the most simply described and easily analyzed.
For this reason, if $\dim(F)>1$, i.e., $\Br(F)\neq 0$, it is logical to formulate
some measure of $F$'s arithmetic complexity in terms of the extent to which all algebras/Brauer classes
are cyclic.

We say $\H^2(F,\mu_n)$ is {\it generated by $\Z/n$-cyclic classes} when
the cup product map 
\[\H^1(F,\mu_n)\otimes_\Z\H^1(F,\Z/n)\longrightarrow\H^2(F,\mu_n)\]
is surjective.
Whether or not $\H^2(F,\mu_n)$ is generated by $\Z/n$-cyclic classes for all $F$ and 
prime-to-$\car(F)$ $n$ is an open problem, posed by Albert in 1936 (\cite[p.126]{Al36}).
It is known
when $F$ contains the $n$-th roots of unity $\mu_n$ or when $n=3$ by Merkurjev-Suslin
\cite[Theorem 16.1, Corollary 16.4]{MS83}, when $n=\ell$ is prime and $[F(\mu_\ell):F]\leq 3$ 
by Merkurjev \cite[Corollary, p.2616]{Mer83}, when $n=5$ by Matzri \cite[Theorem 4.5]{Mat08},
and when $n=\ell$ is prime and every class of {\it index} $\ell$ is cyclic,
by Merkurjev \cite[Theorem 2]{Mer83}.

When $\H^2(F,\mu_n)$ is generated by $\Z/n$-cyclic classes we define
the {\it $\Z/n$-length} $\nL(F)$ to be the smallest $N\in\N\cup\{\infty\}$
such that any class in $\H^2(F,\mu_n)$ can be expressed as a sum of
$N$ $\Z/n$-cyclic classes (see Tignol's paper \cite{Tig84} for background and results on this
invariant, whose origin would seem to go back to Teichm\"uller in \cite{Tei37}).
We define the {\it $n$-Brauer dimension} $\nBrdim(F)$
to be the smallest number $M\in\N\cup\{\infty\}$ such that any class in $\H^2(F,\mu_n)$
has index dividing $n^M$
(see e.g. \cite{HHK09}).
It is easy to see that $\nBrdim(F)\leq\nL(F)$.
Although finite $\Z/n$-length implies finite $n$-Brauer dimension, it is not known
whether finite $n$-Brauer dimension implies finite $\Z/n$-length, or even if finite
$n$-Brauer dimension implies $\H^2(F,\mu_n)$ is generated by $\Z/n$-cyclic classes.

Let $F$ be the function field of a smooth $p$-adic curve.
Saltman showed in \cite[Theorem 3.4]{Sa97} that $\nBrdim(F)=2$ for any prime-to-$p$ number $n$,
and that for $n=\ell\neq p$ a prime, every class of index $n$ is cyclic.
It follows that
${}_\ell\Br(F)$ is generated by $\Z/\ell$-cyclic classes by \cite[Theorem 2]{Mer83}.
Using results of \cite{Sa07},
Suresh showed that ${}_\ell\Length(F)=2$ for a prime $n=\ell\neq p$ assuming $\mu_\ell\subset F$
in \cite{Sur10}, and two of the authors of the present paper proved ${}_\ell\Length(F)=2$
in general in \cite{BT11b}. 

The main result of this paper is that in this situation,
where $F$ is the function field of a $p$-adic curve and $n$ is any number prime to $p$,
${}_n\H^2(F,\mu_n)$ is generated by $\Z/n$-cyclic classes, and $\nL(F)=2$, matching the $n$-Brauer dimension.
The case when $n=\ell\neq p$ is an odd prime gives us a new proof of the result in \cite{BT11b},
without using Saltman's machinery.
As a corollary we see that any $F$-division algebra of period $n$ and index $n^2$
is decomposable, and we obtain a new proof of the main result in \cite{Sa97}
that $\nBrdim(F)=2$.  
In \cite{BMT} the authors constructed indecomposable division algebras of prime power order $\ell^e$ 
($\ell\neq p$)
with index strictly less than $\ell^{2e}$ over function fields of smooth $p$-adic curves.  
Therefore, the result in this paper completes the investigation into which index/period 
combinations admit indecomposable division algebras over function fields of $p$-adic curves.

We prove more generally that if $F$ is the function field of a smooth curve over a complete
discretely valued field $K=(K,v)$ such that $\k(v)$ does not have a special case with respect
to $n$, and $\alpha\in\H^2(F,\mu_n)$, then there exists a $\Z/n$-cyclic
class $\gamma$ and a $\Z/n$-cyclic field extension $L/F$ such that $(\alpha-\gamma)_L$
is unramified.
We also show that for any class $\alpha$ in $\H^2(F,\mu_n)$ there exists
a regular (2-dimensional, projective, flat) model $X/\O_v$ with underlying reduced closed fiber $C$,
with respect to which $\alpha$
may be expressed as a sum of two cyclic classes via characters defined
over $\k(C)$, and a class arising from $\H^2(\k(C),\mu_n)$.
Finally we show that if the function field $\bar F$ of any smooth curve over the residue field
$\k(v)$ satisfies a certain Grunwald-Wang type hypothesis for algebras, then $\H^2(\bar F,\mu_n)$
is generated by $\Z/n$-cyclic classes implies the same for $\H^2(F,\mu_n)$, and if
$\nL(\bar F)$ is finite then
$\nL(F)\leq 2(\nL(\bar F)+1)$.

\section{Background and Conventions}

\Paragraph{\bf General Conventions.}
Let $S$ be an excellent scheme, $n$ a number that is invertible on $S$,
and $\Lambda=(\Z/n)(i)$ the \'etale sheaf $\Z/n$ twisted by an integer $i$.
We write $\H^q(S,\Lambda)$ for the \'etale cohomology group,
and if $\Lambda$ is arbitrary and fixed (or doesn't matter) 
we write $\H^q(S)$ instead of $\H^q(S,\Lambda)$,
and $\H^q(S,r)$ in place of $\H^q(S,\Lambda(r))$.
If $S=\Spec A$ for a ring $A$, then we 
write $\H^q(A)$ and $\H^q(A,r)$.
If $v$ is a valuation on a field $F$, we write $\k(v)$ for
the residue field of the valuation ring $\O_v$, and $F_v$ for the completion of $F$
at $v$.  If $v$ arises from a prime divisor $D$ on $S$, we write $v=v_D$, $\k(D)$,
and $F_D$.
If $T$ is an integral closed subscheme of $S$ we write $\k(T)$ for its function field. 
If $T\to S$ is a morphism of schemes, then
the restriction $\res_{S|T}:\H^q(S)\to\H^q(T)$ is defined,
and we write $\beta|_T$ or $\beta_T$ for $\res_{S|T}(\beta)$.

If $S$ is noetherian,
let $\Div S$ denote the group of Cartier divisors on $S$, and let $\cZ^1(S)$ denote the
group of cycles of codimension one.  Then we have a group homomorphism
\[\Div S\longrightarrow\cZ^1(S)\]
that preserves support on effective divisors (see \cite[Exercise 7.2.4]{Liu}).
If $S$ is regular the map is an isomorphism by \cite[7.2.16]{Liu}.
In this case we identify a divisor $D$ with its associated cycle, 
and if $D$ is effective we identify it with the associated closed subscheme
(see \cite[Section 21.7]{EGAIV:d}).
If $S$ is regular with function ring $F$ then we write $\div(f)$
for the principal divisor determined by $f\in F^*$, that is,
$\div(f)=\sum v_D(f)D$,
where the (finite) sum is over prime divisors on $S$.

\Paragraph{\bf Residues and Ramification.}\label{residues}
All valuations will be discrete of rank one.
If $F=(F,v)$ is a discretely valued field, $F_v$ is the completion of $F$ at $v$,
$n$ is prime to $\car(\k(v))$, and $\Lambda=(\Z/n)(i)$ for some $i$,
then there is a commutative diagram
\[\xymatrix{
0\ar[r]&\H^q(\O_v)\ar[r]\ar[d]&\H^q(F)\ar[r]\ar[d]\ar@{-->}[dr]^-{\partial_v}&\H^{q+1}_z(\O_v)\ar[r]\ar[d]^-\wr&0\\
0\ar[r]&\H^q(\k(v))\ar[r]&\H^q(F_v)\ar[r]&\H^{q-1}(\k(v),-1)\ar[r]&0}\]
The vertical arrow on the right is an isomorphism by cohomological purity (see \cite[Section 3.3]{CT95}).
We call the top row and any (long or short) exact sequence of this type a {\it localization sequence} (see \cite[III.1.25]{M}).
The top row is short exact because $\O_v$ is a discrete valuation ring (see \cite[Section 3.6]{CT95})
We call the bottom row a {\it Witt exact sequence}.
The bottom surjection is split (non canonically) by the map $\omega\mapsto(\pi_v)\cdot\omega$,
where $(\pi_v)\in\H^1(F_v,\mu_n)$ is the Kummer element
determined by a choice of uniformizer $\pi_v$ for $F_v$.
We call the resulting direct sum decomposition $\H^q(F_v)\isom\H^q(\k(v))\oplus\H^{q-1}(\k(v),-1)$
a {\it Witt decomposition}.
We call the map
\begin{equation}\label{residuemap}
\partial_v:\H^q(F)\longrightarrow\H^{q-1}(\k(v),-1)
\end{equation}
the {\it residue map at} $v$.
We will also call any map factoring through $\partial_v$ a residue map.
We say $\alpha\in\H^q(F)$ is {\it unramified at $v$}
if $\partial_v(\alpha)=0$, and {\it ramified at $v$} if $\partial_v(\alpha)\neq 0$.
If $\alpha$ is unramified at $v$ then $\alpha$ comes from $\H^q(\O_v)$,
and we say $\alpha$ is {\it defined} at $\k(v)$.
If $\alpha$ is defined at $\k(v)$ then it has a {\it value} 
\[\alpha(v)=\res_{\O_v|\k(v)}(\alpha)\;\in\H^q(\k(v))\]
Since the ring homomorphism $\O_v\to\k(v)$ factors through the complete discrete valuation
ring $\O_{F_v}\subset F_v$ we have
the alternative description $\alpha(v)=\res_{F|F_v}(\alpha)$, using the Witt sequence
to identify $\H^q(\k(v))$ with the subgroup $\H^q(\O_{F_v})\leq\H^q(F_v)$.
Note $\alpha_{F_v}=0$ if and only if
$\partial_v(\alpha)=0$ and $\alpha(v)=0$ by the Witt sequence.
If $v$ arises from a prime divisor $D$ on an integral scheme
with function field $F$ we generally substitute $D$ for $v$,
and write $\partial_D$
and $\alpha(D)$ in place of $\partial_v$ and $\alpha(v)$.

Each $f\in F^*$ defines a Kummer element $(f)\in\H^1(F,\mu_n)$,
and $\partial_v(f)=v(f)\pmod n\in\H^0(\k(v),\Z/n)$.
If $\chi\in\H^q(F,-1)$ we write $(f)\cdot\chi\in\H^q(F)$ for the cup product.  
Then 
\begin{equation}\label{cupresidue}
\partial_v((f)\cdot\chi)=\left[v(f)\cdot\chi-(f)\cdot\partial_v(\chi)+(-1)\cdot v(f)\cdot\partial_v(\chi)\right]_{F_v}
\end{equation}
See \cite[II.7.12, p.18]{GMS} for the general cup product formula.

Each $\theta\in\H^1(F,\Z/n)$ determines a cyclic field extension, which we denote by $F(\theta)/F$.
Let $w=w_1,\dots,w_g$ be the extensions of $v$ to $F(\theta)$.
Then we have $g$ $\k(v)$-isomorphic residue extensions $\k(w)/\k(v)$ of degree $f$,
a ramification index $e=|\partial_v(\theta)|$ such that $|\theta|=[F(\theta):F]=efg$, and a 
commutative diagram
\[
\xymatrix{
&\H^q(F(\theta))\ar[r]^-{\partial_w}&\H^{q-1}(\k(w),-1)\\
&\H^q(F)\ar[u]^\res\ar[r]_-{\partial_v}&\H^{q-1}(\k(v),-1)\ar[u]_{e\cdot\res}\\
}
\]
If $\theta_{F_v}$ has (Witt) decomposition $\theta_{F_v}=\theta^\circ+(\pi_v)\cdot\partial_v(\theta)$
then $\k(w)=\k(v)(e\cdot\theta^\circ)$.

\Paragraph{\bf Ramification Divisor.}\label{ramdivisor}
If $S$ is a noetherian scheme, $n$ is invertible on $S$, $y\in S$ is a point, and $\alpha\in\H^q(S)$,
Then for each discrete valuation $v$ on $\k(y)$
$\alpha$ has a residue 
$$
\partial_{y,v}(\alpha)\df\partial_v(\alpha_{\k(y)})\in\H^{q-1}(\k(v),-1)
$$
We write $\partial_v(\alpha)$ if $y$ is implied from the context.

More generally, following Kato in \cite{Ka86}, suppose $S$ is a noetherian excellent scheme,
$n$ is invertible on $S$, $y\in S$ is a point
of dimension $d\geq 1$, and $z\in S$ is a point of dimension $d-1$.
We write $y\in S_{(d)}$ and $z\in S_{(d-1)}$.
Let $Y$ be the normalization of $\overline{\{y\}}$ in $\k(y)$.  
Then $Y\to\overline{\{y\}}$ is finite since $S$ is excellent, 
and $\partial_{y,z}$ is defined to be zero if $z\not\in\overline{\{y\}}$, otherwise
\[\partial_{y,z}\;\df\;\sum_{v|z}\cor_{\k(v)|\k(z)}\cdot\partial_{y,v}:\H^q(\k(y))\rightarrow\H^{q-1}(\k(z),-1)\]
where $v$ ranges over the (finitely many) points of $Y$ lying over $z$, and $\partial_{y,v}$
is as in \eqref{residuemap}.
Now if $\alpha\in\H^q(S)$ we have a restriction map $\H^q(S)\to\H^q(\k(y))$,
and we write \[\partial_{y,z}(\alpha)\;\df\;\partial_{y,z}(\alpha_{\k(y)})\in\H^{q-1}(\k(z),-1)\]

Suppose $\alpha\in\H^q(S)$, and $S\to T$ is a birational morphism of 
reduced noetherian excellent schemes (see \cite[I.2.2.9]{EGAI}).
Then $S$ and $T$ have the same generic points (for example, $T$ could be an integral scheme
and $S$ could be its generic point).
Assume every irreducible component of $T$ has dimension $d\geq 1$.
For each $z\in T_{(d-1)}$ write $\partial_z$ for $\bigoplus_{y\in T_{(d)}}\partial_{y,z}$
so that we have a map 
\begin{equation}\label{partialz}
\partial_z:\bigoplus_{y\in T_{(d)}}:\H^q(\k(y))\longrightarrow\H^{q-1}(\k(z),-1)
\end{equation}
The {\it ramification locus of $\alpha$ on $T$} is the (codimension-one) cycle 
\[D_\alpha=\sum_{\underset{\mbox{\scriptsize $\partial_z(\alpha)\neq 0$}}{z\in T_{(d-1)}}}
\overline{\{z\}}\]
We will also use the notation $\ramloc(\alpha)$ in place of $D_\alpha$.
If $T$ is a regular integral scheme then $D_\alpha$ is an effective Cartier divisor, and we call it
the {\it ramification divisor of $\alpha$}.
If $\alpha=(\alpha_x)\in\bigoplus\H^q(\k(x))$ we define the {\it support} of $\alpha$
to be the cycle \[\Supp(\alpha)=\sum_{\alpha_x\neq 0}\overline{\{x\}}\]
We will generally write $\alpha_D$ instead of $\alpha_x$, where $D=\overline{\{x\}}$.

\Paragraph{\bf Kato Complex.}\label{katocomplex}
Let $S$ be an excellent scheme,
$n$ a number invertible on $S$, and $\Lambda=(\Z/n)(i)$ for some $i$.
Suppose $S$ has points of dimension $d\geq 1$.
Identify $\H^{q-1}(\k(z),-1)$ with its image in $\bigoplus_{S_{(d-1)}}\H^{q-1}(\k(x),-1)$,
and define
\[\partial_q=\sum_{z\in S_{(d-1)}}\partial_z:\bigoplus_{x\in S_{(d)}}\H^q(\k(x))\longrightarrow\bigoplus_{x\in S_{(d-1)}}\H^{q-1}(\k(x),-1)\]
where $\partial_z$ is from \eqref{partialz}.
Kato proved in \cite{Ka86} that $\partial_{q-1}\circ\partial_q=0$.
Thus if $S$ is an integral excellent scheme of dimension 2 with function field $F$,
the maps $\partial_q$ define a complex (see \cite{Ka86})
\[\xymatrix{
C_n^{-1}(S):\qquad
\H^2(F)\ar[r]^-{\partial_2}&\bigoplus_{x\in S_{(1)}}\H^1(\k(x),-1)\ar[r]^-{\partial_1}
&\bigoplus_{x\in S_{(0)}}\H^0(\k(x),-2)}\]
Note in particular that if $\alpha\in\H^2(F)$ and $D_\alpha$ contains exactly two
components $D_1=\overline{\{y_1\}}$ and $D_2=\overline{\{y_2\}}$ 
meeting at the closed point $z\in S_{(0)}$, then 
\begin{equation}\label{sumtozero}
\partial_{y_1,z}(\partial_{D_1}(\alpha))+\partial_{y_2,z}(\partial_{D_2}(\alpha))=0
\end{equation}
Since we have a complex we have homology groups; we set 
\[\h_1(C_n^{-1}(S))\;\df\;\ker(\partial_1)/\im(\partial_2)\]

\Paragraph{\bf $\Z/n$-Cyclic Classes.}\label{cycliclength}
If $F$ is any ring we say a class $\alpha\in\H^2(F,\mu_n)$ is {\it $\Z/n$-cyclic} if it is of the form
$(f)\cdot\chi$ for some $(f)\in\H^1(F,\mu_n)$ and $\chi\in\H^1(F,\Z/n)$, and we say 
$\H^2(F,\mu_n)$ is {\it generated by $\Z/n$-cyclic classes} 
if the cup product map
\[\H^1(F,\mu_n)\otimes_\Z\H^1(F,\Z/n)\longrightarrow\H^2(F,\mu_n)\]
is surjective.  

\Definition{}\label{Znlength}
\begin{enumerate}[(a)]
\item
The {\it $\Z/n$-length $\nL(\alpha/F)$ of $\alpha\in\H^2(F,\mu_n)$} is the smallest $M\in\N\cup\{\infty\}$
such that $\alpha$ may be expressed as a sum of $M$ $\Z/n$-cyclic classers.
\item
If $\H^2(F,\mu_n)$ is generated by $\Z/n$-cyclic classes, 
the {\it $\Z/n$-length $\nL(F)$ of $F$} is the smallest $N\in\N\cup\{\infty\}$
such that every class in $\H^2(F,\mu_n)$
may be expressed as a sum of $N$ $\Z/n$-cyclic classes.
Thus $\nL(F)=\sup_\alpha\{\nL(\alpha/F)\}$.
\end{enumerate}

\Paragraph{\bf Basic Setup.}\label{setup}
We use \cite{Liu} as a basic reference.
In all that follows $R$ will be a complete discrete valuation ring with fraction
field $K$ and residue field $k$, $n$ will be a number invertible in $R$, 
$F$ will be the function field of a smooth projective curve over $K$,
and $X/R$ will be a regular (projective, flat) relative curve with function field $F$.
Thus $X$ is 2-dimensional, all of its closed points $z$ have codimension 2 (\cite[8.3.4]{Liu}), 
and the corresponding local rings $\O_{X,z}$ are factorial.
The closed fiber $X_0=X\otimes_R k$ is a connected projective curve over $k$;
we write $C=X_{0,\red}$ for its reduced subscheme, 
$C_1,\dots,C_m$ for the irreducible components of $C$, and $\k(C)=\prod_i\k(C_i)$.
We assume throughout that each irreducible component $C_i$ of $C$ is regular, and that at most
two components meet at a given point, a situation that can always be achieved by blowing up 
(see \cite[9.2.26]{Liu}).
We let $\S$ denote the set of singular points of $C$.
and let $\O_{C,\S}$
denote the localization $\O_{C,\S}=\varinjlim_{\S\subset U}\O_C(U)$.
If $z$ is a closed point of $X$ then $z$ lies on $C$,
and if $z\in C_i$ we write $K_{i,z}=\Frac(\O_{C_i,z}^\h)$, 
where the superscript ``$\h$'' denotes {\it henselization}. 
If $z\in\S$ is on $C_i\cap C_j$ we let $K_z=\Frac(\O_{C,z}^\h)=K_{i,z}\times K_{j,z}$.

Since exactly two irreducible components of $C$ meet at any $z\in\S$
the {\it dual graph} $\G_C$ is defined, and consists of a vertex for each irreducible
component of $C$ and an edge for each singular point, 
such that an edge and a vertex are incident when the corresponding singular point
lies on the corresponding irreducible component (\cite[2.23]{Sai85}, see also \cite[10.1.48]{Liu}). 
The (first) {\it Betti number} $\b_1=\rk(\H_1(\G_C,\Z))$ is the sum $N+E-V$,  
where $V,E$ and $N$ are the numbers of vertices, edges, and connected components of 
$\G_C$.

\Paragraph{\bf Distinguished Divisors.}\label{distinguished}
We will want to compute the residues of elements of $\H^q(F)$ lifted from $\H^q(\k(C))$,
and to do it we need 
a mechanism for lifting 0-cycles on $C$, which are finite sets of closed points, 
to divisors on $X$.
See \cite[Section 8.3]{Liu} for background on the following definitions.

Assume the setup of \eqref{setup}.
We say a divisor $D\in\Div X$ is {\it horizontal} if each irreducible component of $\Supp(D)$
maps surjectively to $\Spec R$.
Since $X/R$ is projective, it follows by Zariski's main theorem
that an effective horizontal divisor is finite over $\Spec R$, and hence,
since $R$ is henselian, that any irreducible
effective horizontal divisor has a single closed point (\cite[I.4.2]{M}).
Since we assume each $C_i$ is regular and exactly two components meet at any $z\in\S$,
for each closed point $z\in C$ there exists a regular effective horizontal prime divisor $D$ on $X$ passing through $z$
that is transverse at $z$ to each irreducible component of $C$ passing through $z$
(\cite[Proposition 2.6]{BT11}).
Let $\D$ denote a preassigned set of these divisors, and call them {\it distinguished prime divisors}.
Let $\D_\S$ denote the subset of those that avoid $\S$,
and let $\D_*$ denote the subset of those that pass through $\S$.
We will abuse notation and say an effective divisor $D$ is ``in $\D$'' if it is
a sum of distinguished prime divisors (with multiplicity one), and if $D$ is any 
divisor, we denote by $\D(D)$ (or $\D_\S(D)$, etc.) the 
set of prime components of $D$ that are in $\D$ (or $\D_\S$, etc.).

Once we have fixed a set $\D$,
$D_z$ will denote the distinguished prime divisor passing through a closed point $z$ of $C$.
Since $D_z$ is finite over $\Spec R$ it is affine, and if $\pi_{D_z}\in\O_{X,z}$ is a local equation
for $D_z$ then $D_z=\Spec\O_{X,z}/(\pi_{D_z})$.
Since $D_z$ is regular $\O_{X,z}/(\pi_{D_z})$
is a complete discrete valuation ring whose fraction field $\k(D_z)$ is complete 
with respect to the discrete valuation defined by $z$.  
Since $D_z$ is transverse to each irreducible component $C_k$ passing through $z$, 
the image in $\k(D_z)$ of any local equation $\pi_{C_k,z}\in\O_{X,z}$ for $C_k$ is a uniformizer for $\k(D_z)$.
Conversely the image in $\k(C_k)$ of $\pi_{D_z}\in\O_{X,z}$ is a uniformizer for the 
discrete valuation on $\k(C_k)$ defined by $z$.

\Lemma\label{picard}{
Assume the setup of \eqref{setup}.
Let $\T$ be any 0-cycle on $C$, and let $\D_\T$
denote the distinguished prime divisors that avoid $\T$.
Then $\Pic X/n$ is generated by (the classes of) a subset of $\D_\T$,
and intersection with $C$ determines an isomorphism
$\Pic X/n\isom\Pic C/n$.
}
\rm

\begin{proof}
(See e.g. \cite{Sa08}.)
We may assume that $\T$ contains $\S$ and includes a point from each irreducible component of $C$.
Let $\mathscr H_\T\subset\Div X$ denote the subgroup generated by horizontal prime Cartier divisors on $X$
that avoid $\T$.
Since $X$ and $C$ are noetherian with no embedded points (by \eqref{setup}),
the maps $D\mapsto\O_X(D)$ and $D_0\mapsto\O_C(D_0)$ induce
isomorphisms $\CaCl(X)\isom\Pic X$ and $\CaCl(C)\isom\Pic C$, by \cite[7.1.19]{Liu}.

We first show $\CaCl(X)$ is generated by elements of $\mathscr H_\T$.
Fix $D\in\Div X$.
Since $\T$ is finite and $X$ is regular, 
the semi-local ring $\O_{X,\T}$ is factorial.
Since $\Spec(\O_{X,\T})\to X$ is a (birational) localization, 
the restriction of $D$ to $\Spec(\O_{X,\T})$ is a Cartier divisor,
and since $\O_{X,\T}$ is factorial, 
$D|_{\O_{X,\T}}$ is principal, hence $D|_{\O_{X,\T}}=\div(f)$ for some $f\in F$.
Since $\div(f)$ matches $D$ on $\Spec\O_{X,\T}$,
the support of $D-\div(f)\in\Div X$ avoids $\T$.
Since $\T$ includes a point from each component of $C$, $\Supp(D-\div(f))$ contains no
irreducible component of $C$, hence
$D-\div(f)$ is in $\mathscr H_\T$, as desired.
Thus any class $[D]\in\CaCl(X)\isom\Pic X$ may be represented by 
an element of $\mathscr H_\T$.

Since the support of elements of $\mathscr H_\T$ contain no generic points of $C$,
the restriction map $\mathscr H_\T\to\Div C$ given by $D\mapsto D|_C$ is well defined
and sends principal divisors to principal divisors by \cite[Lemma 7.1.29]{Liu}.
Since $\mathscr H_\T$ generates $\CaCl(X)$ it defines a map $\CaCl(X)\to\CaCl(C)$.
This map is compatible with the natural map $\Pic X\to\Pic C$, i.e.,
if $i:C\to X$ is the closed immersion, then
$i^*\O_X(D)=\O_C(D|_C)$ for $D\in \mathscr H_\T$.  This is easily seen
by looking at trivializing neighborhoods (see also \cite[Lemma 7.1.29]{Liu}).

We claim $\D_\T\subset \mathscr H_\T$ generates a subgroup
of $\Pic X$ that maps onto $\Pic C$.
Since $X/R$ is projective and $R$ is henselian, for each effective (Cartier) divisor $D_0\in\Div C$
there exists a divisor $D\in\Div X$ such that $D|_C=D_0$,
by \cite[Proposition 21.9.11]{EGAIV:d} and the proof of
\cite[Corollary 21.9.12]{EGAIV:d}.
As before, there exists an $f\in F^*$ such that $D-\div(f)$ is in $\mathscr H_\T$.
Since $D$ and $D-\div(f)$ are horizontal,
$f$ is regular at every generic point of $C$, hence we have an image $f_0\in\k(C)^*$,
hence $(D-\div(f))|_C=D_0-\div(f_0)$ is a divisor 
that is equivalent to $D_0$.
Since $\Supp(D-\div(f))$ avoids $\T$, so does $\Supp(D_0-\div(f_0))$ by \cite[7.1.29]{Liu}.
Thus $[D_0]$ may be represented by a Cartier divisor $D_0'$
that avoids $\T$, and since $\T$ includes $\S$ we may lift $D_0'$ to a divisor $D'\in\D_\T$
by \cite[Proposition 2.6]{BT11} (see also \cite[8.3.35(b)]{Liu}).

We show $\Pic X\to\Pic C$ is bijective modulo $n$.
By the Kummer sequence we have a diagram
\begin{displaymath}
\xymatrix{
0 \ar[r] & \Pic X/n \ar[d]\ar[r] & \H^2(X,\mu_n)\ar[d]\ar[r] & {}_n\H^2(X,\Gm)\ar[d]\ar[r] & 0\\
0 \ar[r] & \Pic C/n \ar[r] & \H^2(C,\mu_n) \ar[r] & {}_n\H^2(C,\Gm)\ar[r] & 0}
\end{displaymath}
Since $R$ is complete, the middle down-arrow is an isomorphism by proper base change,
so $\Pic X/n\to\Pic C/n$ is injective, hence bijective.
Finally, since the subgroup of $\Pic X$ generated by $\D_\T$ maps onto $\Pic C$,
this shows $\D_\T$ generates $\Pic X/n$.

\end{proof}

\Remark\label{picardtree}
If $F=K(T)$ is a rational function field over $K$ then $X$ may be obtained
from $\P^1_R$ as a series of blowups at closed points, and then $C$ is a tree
with irreducible components $C_i=\P^1_{\k(z_i)}$ for various fields $\k(z_i)$ finite over $k$,
by \cite[Theorem 1.19]{Liu}.
Then by \cite[Proposition 2.18]{Liu} 
$\Pic X\isom\Z^m$, where $m$ is the number of irreducible components of $C$,
and we may generate $\Pic X$ by choosing horizontal divisors $D_i\in\D_\S$, one
passing through each $C_i$.
Thus in this case $\Pic X/n\isom\Pic C/n$ is finite, with a single pre-assigned
horizontal generating divisor for each irreducible component of $C$.

\Corollary\label{picardcor}
Assume the setup of \eqref{setup}, and suppose given $E\in\Div X$.
Then there exists an element $f\in F$ such that $\div(f)=D_{(f)}=E+H\pmod n$, for some $H\in\D_\S$
(possibly with multiplicities) that avoids any finite preassigned set of closed points.  
\rm

\begin{proof}
This is immediate from Lemma~\ref{picard}.

\end{proof}

\section{Gluing Theory}

Assume the basic setup \eqref{setup}, with $\Lambda=(\Z/n)(i)$. 
In \cite{BT11} the authors proved a ``gluing lemma'' \cite[Lemma 4.2]{BT11}
that defined a compatibility criterion for the entries
of a tuple in the product $\H^q(\k(C))=\prod_i\H^q(\k(C_i))$ to
be in the image of the natural map $\H^q(\O_{C,\S})\to\H^q(\k(C))$. 
They then constructed (in \cite[Theorem 4.9]{BT11})
a homomorphism $\H^q(\O_{C,\S})\to\H^q(F)$ for all $q\geq 0$,
providing a means of constructing elements of $\H^q(F)$ with known properties,
leading to the construction of noncrossed product $F$-division algebras
and indecomposable $F$-division algebras of unequal period and index (see also \cite{BMT}).
Our goal in this section is to improve the gluing lemma \cite[Lemma 4.2]{BT11}
by expanding the compatibility criterion to a larger subgroup of $\H^q(\k(C),\Lambda)$
(denoted by $\Gamma_{\D_*}^q(\k(C),\Lambda)$) 
that includes tuples whose entries are ramified at the singular points $\S$ of $C$.
This is accomplished in Theorem~\ref{newglue} below.
The access to tuples that ramify at $\S$ allows us
to ``reach'' any class of $\H^2(F,\mu_n)$ over some model $X$ using concrete constructions
over $\k(C)$, and as a result we are able to
bound the $\Z/n$-length of $\H^2(F,\mu_n)$.
These applications are treated in sections 3, 4, and 5.

\Paragraph{\bf Gluing Framework.}\label{gluingframework}
Assume the setup of \eqref{setup}.
Let $D$ be a prime divisor in $\D$ and let $U=X-D$ be the open subscheme.
By \cite[III.1.25]{M} we have a long exact localization sequence
\[\xymatrix{
\cdots\ar[r]&\H^q(X)\ar[r]&\H^q(U)\ar[r]&\H^{q+1}_D(X)\ar[r]&\H^{q+1}(X)\ar[r]&\cdots}\]

Since $X$ and $D$ are regular we have a natural isomorphism $\H^{q+1}_D(X)\isom\H^{q-1}(D,-1)$
by Gabber's absolute purity theorem \cite{Fuj02}, and the induced map
$\H^q(U)\to\H^{q-1}(D,-1)$ is the residue map $\partial_D$.

Let $V=X-\D$.  This may not be a scheme, but we use the notation heuristically,
and write \[\H^q(V)\;\df\;\varinjlim\H^q(U)\] where the limit is over open $U$
such that $X-U\in\D$.  
Taking the limit over the localization sequences and setting $\partial=\bigoplus_\D\partial_D$ yields
\[\xymatrix{
\cdots\ar[r]&\H^q(X)\ar[r]&\H^q(V)\ar[r]^-{\partial}&\bigoplus_\D\H^{q-1}(D,-1)\ar[r]&
\H^{q+1}(X)\ar[r]&\cdots}\]

Since $C/k$ is projective it is separated (over $\Z$), 
hence if $U_0\subset C$ is an affine open subset of $C$ then $U_0\to C$ is an affine map by \cite[3.3.6]{Liu}.
The inverse limit $\varprojlim U_0$ over dense affine open subschemes of $C$ is then a scheme
by \cite[8.2.3]{EGAIV:c}, and this scheme is $\Spec\k(C)$ by \cite[Exercise 5.1.17(c)]{Liu}.
Therefore $\H^q(\k(C),\Lambda)=\varinjlim\H^q(U_0,\Lambda)$ by \cite[III.1.16]{M}, where the limit
is over dense open $U_0\subset C$.

Since $X$ is regular at each point, the closed fiber $X_0$ has no embedded points.
For if $p=u\pi_1^{e_1}\pi_2^{e_2}$ is a factorization of the uniformizer of $R$
in the (factorial) ring $\O_{X,z}$
then $\Ass(\O_{X_0,z}=\O_{X,z}/(u\pi_1^{e_1}\pi_2^{e_2}))$ is just 
the set of minimal primes $\{(\pi_1),(\pi_2)\}$.
It follows that $\k(C)$ equals $\k(X_0)_\red$ by \cite[Lemma 4.8]{BT11}.
Now since $C=X_{0,\red}$ and $\k(C)=\k(X_0)_\red$ 
we have isomorphisms 
\[\H^q_z(X_0)=\H^q_z(C)\,,\quad\H^q(X_0)=\H^q(C)\,,\quad\text{and}\quad\H^q(\k(X_0))=\H^q(\k(C))\]
by \cite[V.2.4(c)]{M} (see also \cite[II.3.11]{M}).
Since $R$ is complete we then have isomorphisms $\H^q(X)\isim\H^q(C)$
and $\H^q(D)\isim\H^q(\k(z))$ for $D=D_z\in\D$ by \cite[VI.2.7]{M}.
Then, because of the functorial nature of the localization sequence \cite[III.1.25]{M}
(and the exactness of the direct limit functor),
we obtain for any $q\geq 0$ and fixed $\Lambda$ a commutative ladder
\begin{equation}\label{maplocseqs}
\xymatrix{
\cdots\ar[r]&\H^q(X)\ar[r]\ar[d]^\wr&\H^q(V)\ar[r]^-{\;\partial\;}\ar[d]&\bigoplus_{\D}\H^{q-1}(D,-1)
\ar[r]\ar[d]&\H^{q+1}(X)\ar[d]^\wr\ar[r]&\cdots\\
\cdots\ar[r]&\H^q(C)\ar[r]&\H^q(\k(C))\ar[r]^-{\;\delta\;}&
\bigoplus_{C_{(0)}}\H^{q+1}_z(C)\ar[r]&\H^{q+1}(C)\ar[r]&\cdots
}
\end{equation}
The goal of gluing is to invert the map $\H^q(V)\to\H^q(\k(C))$ (on its image).
If $z\in C$ is a closed point then we have the localization sequence
\[\xymatrix{\cdots\ar[r]&\H^q(C)\ar[r]&\H^q(C-\{z\})\ar[r]&
\H^{q+1}_z(C)\ar[r]&\H^{q+1}(C)\ar[r]&\cdots}\]
Since
$\H^{q+1}_z(C)=\H^{q+1}_z(\O_{C,z}^\h)$ by excision \cite[III.1.28]{M} we obtain
an ``excised'' localization sequence
\begin{equation}\label{exs}
\xymatrix{\cdots\ar[r]&\H^q(\O_{C,z}^\h)\ar[r]&\H^q(K_z)\ar[r]^-{\hat\delta_z}&
\H^{q+1}_z(\O_{C,z}^\h)\ar[r]&\H^{q+1}(\O_{C,z}^\h)\ar[r]&\cdots}
\end{equation}
If $z$ is a normal point then $\O_{C,z}^\h$ is a discrete valuation ring with henselian valued
fraction field $K_z=\Frac\O_{C,z}^\h$, and
then $\H^{q+1}_z(C)\isom\H^{q-1}(\k(z),-1)$ by Witt's theorem.
Thus if $C$ is nonsingular we can replace $\H^{q+1}_z(C)$
with $\H^{q-1}(\k(z),-1)$ in the bottom row of \eqref{maplocseqs},
and use the proper base change isomorphism \[H^{q-1}(D,-1)\isom\H^{q-1}(\k(z),-1)\]
and the 5-Lemma to invert the arrow $\H^q(V)\to\H^q(\k(C))$.
But in general $F$ may not admit a model $X$
with an irreducible closed fiber, 
and then $C$ has singularities (i.e., $\S\neq\varnothing$)
and this strategy fails.
Moreover, even if $F$ admits a model $X$ with a nonsingular closed fiber
(e.g. if $F=\Q_p(t)$), we will
have reason to blow up $X$,
thereby introducing singularities to $C$.
When $C$ has singularities it turns out the map $\H^q(V)\to\H^q(\k(C))$ 
is neither surjective nor injective in general (see below),
and this constitutes an obstruction to lifting.
To solve this problem we will characterize the image
of $\H^q(V)$ in $\H^q(\k(C))$ (see Definition~\ref{glueatz} and the proof of Theorem~\ref{newglue}), 
and show the kernel is inconsequential (see Lemma~\ref{cslemma}).
We start with a technical lemma. 

\Lemma\label{keylemma}
Assume \eqref{setup}.
Fix $z\in C_i$ and $D=D_z\in\D$. 
Let $\pi_D\in\O_{X,z}$ be a local equation for $D$,
and let $\bar\pi_D$ be the image of $\pi_D$ in $K_{i,z}=\Frac\O_{C_i,z}^\h$.
Then there is a commutative diagram of compatibly-split exact sequences
$$
\xymatrix{
0\ar[r]&\H^q(\k(D))\ar[r]&\H^q(F_D)\ar[r]&\H^{q-1}(\k(D),-1)
\ar[r]&0\\
0\ar[r]&\H^q(\k(z))\ar[r]\ar[u]^{\inf}\ar@{=}[d]&\H^q(\O_{X,z}^\h[1/\pi_D])
\ar[r]\ar@{>->}[u]\ar[d]^\wr&\H^{q-1}(\k(z),-1)\ar[r]\ar[u]^\inf\ar@{=}[d]&0\\
0\ar[r]&\H^q(\k(z))\ar[r]&\H^q(K_{i,z})\ar[r]&\H^{q-1}(\k(z),-1)
\ar[r]&0
}
$$
The top two splittings are defined by $\omega\mapsto(\pi_D)\cdot\omega$,
the bottom by $\omega\mapsto(\bar\pi_D)\cdot\omega$.
\rm

\begin{proof}
The top and bottom sequences are well known to be split exact (they are ``Witt sequences'')
as indicated.
Since $D$ is in $\D$ the element $\pi_D\in\O_{X,z}$ remains prime in $\O_{X,z}^\h$,
and letting $D$ denote the corresponding closed subset of $\Spec\O_{X,z}^\h$
we compute $\Spec\O_{X,z}^\h-D=\Spec\O_{X,z}^\h[1/\pi_D]$.
The middle sequence is now derived from the localization sequence
\[
\xymatrix{
\cdots\ar[r]&\H^q(\O_{X,z}^\h)\ar[r]&\H^q(\Spec\O_{X,z}^\h[1/\pi_D])\ar[r]&\H^{q+1}_D(\O_{X,z}^\h)\ar[r]&\cdots
}\]
Since $\Spec \O_{X,z}^\h$ and $D$ are regular we compute $\H^{q+1}_D(\O_{X,z}^\h)=\H^{q-1}(D,-1)$ by 
Gabber's absolute purity theorem.
Therefore we have an exact sequence
\begin{equation}\label{xs}
\xymatrix{\cdots\ar[r]&\H^q(\O_{X,z}^\h)\ar[r]&\H^q(\O_{X,z}^\h[1/\pi_D])\ar[r]^-{\partial_D}
&\H^{q-1}(D,-1)\ar[r]&\cdots}
\end{equation}
An element $\omega\in\H^{q-1}(D,-1)$ may be defined over $\O_{X,z}^\h[1/\pi_D]$ using the isomorphisms
$$
\H^{q-1}(D,-1)\longrightarrow\H^{q-1}(\k(z),-1)\longleftarrow\H^{q-1}(\O_{X,z}^\h,-1)
$$
which follow from \cite[Remark III.3.11(a)]{M}.
Therefore since $\partial_D((\pi_D)\cdot\omega)=\omega$ by \eqref{cupresidue}, 
the map $\partial_D:\H^q(\O_{X,z}^\h[1/\pi_D])\to\H^{q-1}(D,-1)$ is surjective, and
is split by the map $\omega\mapsto(\pi_D)\cdot\omega$.
It follows from the exactness of \eqref{xs} that the map $\H^q(\O_{X,z}^\h)\to\H^q(\O_{X,z}^\h[1/\pi_D])$ 
in \eqref{xs} is injective.
Identifying $\H^q(\O_{X,z}^\h)$ with $\H^q(\k(z))$ gives the middle split exact sequence.

Since $\O_{C_i,z}^\h$ is a quotient of $\O_{X,z}^\h$,
there is a natural map $\O_{X,z}^\h[1/\pi_D]\to K_{i,z}$, and
by the functoriality of the localization sequence
we have a commutative diagram
$$
\xymatrix{
0\ar[r]&\H^q(\O_{X,z}^\h)\ar[r]\ar[d]^\wr&\H^q(\O_{X,z}^\h[1/\pi_D])\ar[r]\ar[d]&
\H^{q-1}(D,-1)\ar[d]^\wr\ar[r]&0\\
0\ar[r]&\H^q(\O_{C_i,z}^\h)\ar[r]&\H^q(K_{i,z})\ar[r]&\H^{q-1}(\k(z),-1)\ar[r]&0
}
$$
The left vertical arrow is an isomorphism 
since $\H^q(\O_{X,z}^\h)\to\H^q(\k(z))$ factors (isomorphically) 
through the isomorphism $\H^q(\O_{C_i,z}^\h)\to\H^q(\k(z))$.
Since $D$ meets $C_i$ transversely, $\bar\pi_D$ is a uniformizer for $K_{i,z}$,
and the Witt decomposition of the bottom with respect to $\bar\pi_D$ is compatible with the
splitting of the top with respect to $\pi_D$.
Thus the middle vertical arrow is an isomorphism.
Since the groups on the ends are all isomorphic to the corresponding groups over $\k(z)$, this gives
the bottom half of the commutative ladder.

We prove there exists a map $\O_{X,z}^\h[1/\pi_D]\to F_D$. 
Suppose then that $Y\to X$ is a pointed \'etale neighborhood of $z$, so there exists a lift $z\to Y$.
Then $D_Y\to D$ is \'etale, hence there exists a section $D\to D_Y\subset Y$ by \cite[I.4.2(d)]{M},
hence $Y$ is a pointed \'etale neighborhood of $\xi_D:=\Spec\k(D)$.
Thus every pointed \'etale neighborhood of $z$ is naturally
a pointed \'etale neighborhood of $\xi_D$, and we have a map $\O_{X,z}^\h\to\O_{X,\xi_D}^\h$.
Since $\O_{F_D}$ is henselian and $\O_{X,\xi_D}\to\O_{F_D}$ is a local homomorphism
we obtain a map $\O_{X,\xi_D}^\h\to\O_{F_D}$ 
by the universal property of the henselization,
hence a map $\O_{X,z}^\h\to\O_{F_D}$,
hence the desired map $\O_{X,z}^\h[1/\pi_D]\to F_D$, since $\pi_D$ is invertible in $F_D$.

Now we have natural commutative diagrams 
$$
\xymatrix{
\H^q(\O_{X,z}^\h)\ar[r]\ar[d]^\wr&\H^q(\O_{F_D})\ar[d]^\wr
\ar@{}[r]&\H^{q-1}(D,-1)\ar[r]\ar[d]^\wr&\H^{q-1}(\k(D),-1)\ar@{=}[d]\\
\H^q(\k(z))\ar[r]^-\inf&\H^q(\k(D))\ar@{}[r]&\H^{q-1}(\k(z),-1)\ar[r]^-\inf
&\H^{q-1}(\k(D),-1)
}
$$
yielding the top commutative ladder in the statement of the lemma.
Finally, since $\pi_D$ is a uniformizer for $F_D$,
a splitting of the Witt sequence for $F_D$ is given by $\omega\mapsto(\pi_D)\cdot\omega$,
compatible with the splitting of the sequence for $\O_{X,z}^\h[1/\pi_D]$.

\end{proof}

\Lemma\label{subgroup}
Assume \eqref{setup}.
Fix $z\in\S$ on $C_i\cap C_j$ and $D=D_z\in\D$ containing $z$.
Let $\pi_D\in\O_{X,z}$ be a local equation for $D$, with images $\pi_j\in K_{i,z}$
and $\pi_i\in K_{j,z}$.
We have a commutative diagram
\[\xymatrix{
&\H^q(\O_{X,z}^\h[1/\pi_D])\ar[r]^-{\partial_D}\ar[d]^g&\H^{q-1}(D,-1)\ar[d]^h\\
&\H^q(K_{i,z})\times\H^q(K_{j,z})\ar[r]^-{\hat\delta_z}&\H^{q+1}_z(C)
}\]
where the top row is \eqref{xs} and the bottom row is \eqref{exs}.
Then 
\begin{enumerate}[\rm a)]
\item
$\im(h)\isom\H^{q-1}(\k(z),-1)$ and $h$ is injective.
\item
$\im(g)=\{(\chi^\circ+(\pi_j)\cdot\omega,\chi^\circ+(\pi_i)\cdot\omega):
\chi^\circ\in\H^q(\k(z)),\omega\in\H^{q-1}(\k(z),-1)\}$
\item
$\im(g)=\hat\delta_z^{-1}(\im(h))$
\end{enumerate}
\rm

\begin{proof}
We identify $\H^{q-1}(D,-1)$ with $\H^{q-1}(\k(z),-1)$ as in the proof of Lemma~\ref{keylemma}
and claim there is a diagram with exact rows
\begin{equation}\label{deltadiagram}
\xymatrix{
\H^q(\k(z))\;\;\ar@{>->}[r]\ar@{=}[d] &\H^q(\O_{X,z}^\h[1/\pi_D])\ar@{>>}[r]^-{\partial_D}\ar[d]^-g
&\H^{q-1}(\k(z),-1)\ar[d]^h\\
\H^q(\k(z))\;\;\ar@{>->}[r]\ar[d]_{(\id,\id)} 
&\H^q(K_{i,z})\times\H^q(K_{j,z})\ar@{>>}[r]^-{\hat\delta_z}\ar@{=}[d]
&\H^{q+1}_z(C)\ar@{>>}[d]^p\\
\H^q(\k(z))\times\H^q(\k(z))\;\;\ar@{>->}[r]&\H^q(K_{i,z})\times\H^q(K_{j,z})
\ar@{>>}[r]^-{(\partial_z,\partial_z)\;\;}&\H^{q-1}(\k(z),-1)\times\H^{q-1}(\k(z),-1)
}
\end{equation}
The top and bottom rows are the middle and bottom rows of Lemma~\ref{keylemma},
with the bottom row the product for $K_{i,z}$ and $K_{j,z}$.
The middle row is derived from the excised localization sequence
\eqref{exs} at the singular closed point $z$.
For since $\O_{C,z}^\h$ is 1-dimensional and $z$ is the intersection of $C_i$ and $C_j$,
the underlying topological space of $\Spec\O_{C,z}^\h$ is $\{\Spec(K_{i,z}),\Spec(K_{j,z}),z\}$,
therefore $\Spec\O_{C,z}^\h-\{z\}=\Spec(K_{i,z}\times K_{j,z})$,
and we obtain the sequence
\[\xymatrix{
\cdots\ar[r]&\H^q(\O_{C,z}^\h)\ar[r]&\H^q(K_{i,z})\times\H^q(K_{j,z})\ar[r]&\H^{q+1}_z(C)\ar[r]&\cdots}\]
Since $\O_{C,z}^\h$ is henselian we have
$\H^q(\O_{C,z}^\h)=\H^q(\k(z))$ by \cite[Remark III.3.11(a)]{M}, as in Lemma~\ref{keylemma}, 
and this gives the middle row.
The map $\H^q(\k(z))\to\H^q(K_{i,z})\times\H^q(K_{j,z})$ is inflation on each factor by the derivation
of the localization sequence,
and it is injective by Witt's theorem.
Therefore $\hat\delta_z$ is surjective.
The vertical maps are of a functorial nature,
with $g$ the product of the restriction maps.
The diagram commutes by the functoriality of the localization sequences,
and the map $p$ is surjective by diagram chase.
This proves the claim.

The composition $p\circ h$ maps $\H^{q-1}(\k(z),-1)$ isomorphically onto the diagonal
subgroup $\{(\omega,\omega)\}\leq\H^{q-1}(\k(z),-1)\times\H^{q-1}(\k(z),-1)$ by Lemma~\ref{keylemma}.
In particular $h$ is injective and
\[\im(h)\isom\H^{q-1}(\k(z),-1)\] 
proving (a). 
The commutativity of the top right square of \eqref{deltadiagram} together with the surjectivity of $\partial_D$ 
shows $\im(g)\subset\hat\delta_z^{-1}(\im(h))$, proving half of (c).
By Lemma~\ref{keylemma} any element of $\H^q(\O_{X,z}^\h[1/\pi_D])$ has the form
$\chi^\circ+(\pi_D)\cdot\omega$ for $\chi^\circ\in\H^q(\k(z))$ and $\omega\in\H^{q-1}(\k(z),-1)$,
and since $\pi_j$ and $\pi_i$ are the images of $\pi_D$ in $K_{i,z}$ and $K_{j,z}$, respectively, 
we compute 
\[g(\chi^\circ+(\pi_D)\cdot\omega)
=(\chi^\circ+(\pi_j)\cdot\omega,\chi^\circ+(\pi_i)\cdot\omega)\in\H^q(K_{i,z})\times\H^q(K_{j,z})\]
Since any element of the form $(\chi^\circ+(\pi_j)\cdot\omega,\chi^\circ+(\pi_i)\cdot\omega)$
is obtainable in this way using Lemma~\ref{keylemma}, this proves (b).  
To prove the rest of (c) it remains to show every element of $\delta_z^{-1}(\im(h))$ has this form.

Suppose $(\chi_{i,z},\chi_{j,z})\in\H^q(K_{i,z})\times\H^q(K_{j,z})$
and $\hat\delta_z((\chi_{i,z},\chi_{j,z}))=\varepsilon$ is in $\im(h)$.
Then $\partial_z(\chi_{i,z})=\partial_z(\chi_{j,z})=\omega$ for some $\omega\in\H^{q-1}(\k(z),-1)$
by the bottom row of \eqref{deltadiagram},
and as above $(\pi_D)\cdot\omega\in\H^q(\O_{X,z}^\h[1/\pi_D])$ maps to $\varepsilon$ via $h\circ\partial_D$.
Therefore since $g((\pi_D)\cdot\omega)=((\pi_j)\cdot\omega,(\pi_i)\cdot\omega)$ we have 
$(\chi_{i,z}-(\pi_j)\cdot\omega,\chi_{j,z}-(\pi_i)\cdot\omega))\in\ker(\hat\delta_z)$.
We conclude $(\chi_{i,z}-(\pi_j)\cdot\omega,\chi_{j,z}-(\pi_i)\cdot\omega))$
is in $\im(\id,\id)$ by the bottom row of \eqref{deltadiagram}, hence it is equal to an element of the
form $(\chi^\circ,\chi^\circ)\in\H^q(K_{i,z})\times\H^q(K_{j,z})$ for some $\chi^\circ\in\H^q(\k(z))$.
Thus $(\chi_{i,z},\chi_{j,z})=(\chi^\circ+(\pi_j)\cdot\omega,\chi^\circ+(\pi_i)\cdot\omega)$,
which shows $\im(g)\supset\delta_z^{-1}(\im(h))$, completing the proof of (c).

\end{proof}

\Remark
The group $\im(g)\leq\H^q(K_{i,z})\times\H^q(K_{j,z})$ consists of pairs
whose entries {\it glue} along $D_z\in\D_*$ at $z\in\S$, in the sense that
they arise from elements defined in a henselian neighborhood of $z$ on $X$,
by Lemma~\ref{subgroup}.
We now extend this notion to find the tuples in $\H^q(\k(C))=\prod_i\H^q(\k(C_i))$ 
whose entries glue in this sense
along $\D_*$ at all points of $\S$.

\Definition\label{glueatz}
Let $\delta:\H^q(\k(C))\to\bigoplus_{C_{(0)}}\H^{q+1}_z(C)$ be the map of \eqref{maplocseqs},
and let $\delta_z$ be the $z$-component.
For any subset $\mathcal F\subset C_{(0)}$
let $\delta_{\mathcal F}=\sum_{z\in\mathcal F}\delta_z$, and
let $\E_{\mathcal F}=\{D_z:z\in\mathcal F\}\subset\D$.
Consider $\delta_{\mathcal F}$ together with the direct sum of diagrams of Lemma~\ref{subgroup}:
\[
\xymatrix{&\bigoplus_{\mathcal F}\H^q(\O_{X,z}^\h[1/\pi_{D_z}])\ar[d]_{\oplus g}\ar[r]^-{\oplus\partial_{D_z}}
&\bigoplus_{\E_{\mathcal F}}\H^{q-1}({D_z},-1)\ar@{>->}[d]^{\oplus h}\\
&\bigoplus_{\mathcal F}\H^q(K_z)\ar[r]^-{\oplus\hat\delta_z}&\bigoplus_{\mathcal F}\H^{q+1}_z(C)\\
&\H^q(\k(C))\ar[u]^-{\oplus\res}\ar[ur]_-{\delta_{\mathcal F}}
}\]
Then set 
\begin{enumerate}
\item
$\Delta_{\E_{\mathcal F}}^{q+1}(C,\Lambda)=\im(\oplus h)\leq \bigoplus_{\mathcal F}\H^{q+1}_z(C,\Lambda)$
\item
$\Gamma_{\E_{\mathcal F}}^q(\k(C),\Lambda)=\delta_{\mathcal F}^{-1}(\Delta_{\E_{\mathcal F}}^{q+1})\leq \H^q(\k(C),\Lambda)$
\end{enumerate}
If $\E_{\mathcal F}=\{D\}$ is a single prime divisor, write $\Delta_D^{q+1}$ and $\Gamma_D^q$ instead.
If $\mathcal F\subset\S$ we say the entries of the tuples in
$\Gamma_{\E_{\mathcal F}}^q(\k(C))\leq\H^q(\k(C))$ {\it glue at $\mathcal F$ along $\E_{\mathcal F}$}.
In particular that the entries $\alpha_{C_i}$ of tuples $\alpha_C=(\alpha_{C_i})\in\Gamma_{\D_*}^q(\k(C))$
glue at $\S$ along $\D_*$.

\Remark\label{delta}
Assume the setup of \eqref{setup}.
\noindent
\begin{enumerate}[\rm a)]
\item
If $z\not\in\S$ then $\Delta_{D_z}^{q+1}(C)=\H_z^{q+1}(C)\isom\H^{q-1}(\k(z))$,
the corresponding map $h:\H^{q-1}(D,-1)\to\H_z^{q+1}(C)$ is an isomorphism,
and $\Gamma_{D_z}^q(\k(C))=\H^q(\k(C))$.
Thus if $\mathscr E$ contains $\D_*$ then
$\Gamma_{\mathscr E}^q=\Gamma_{\D_*}^q=\bigcap_{{D_z}\in\D_*}\Gamma_{D_z}^q$.
\item
For each $z\in\S$ on $C_i\cap C_j$ suppose $\pi_{D_z}\in\O_{X,z}$ has
images $(\pi_j)\in\H^1(K_{i,z},\mu_n)$ and $(\pi_i)\in\H^1(K_{j,z},\mu_n)$.
Then by Lemma~\ref{subgroup} the elements of $\Gamma_{\D_*}^q(\k(C),\Lambda)$ are precisely the tuples 
$\chi_C=(\alpha_{C_k})\in\H^q(\k(C),\Lambda)$ that have ``equal'' Witt decompositions at each $z\in\S$
with respect to $(\pi_j)$ and $(\pi_i)$:
\begin{align*}
\chi_{C_i,z}&=\chi_z^\circ+(\pi_j)\cdot\omega\\
\chi_{C_j,z}&=\chi_z^\circ+(\pi_i)\cdot\omega
\end{align*}
for some $\chi_z^\circ\in\H^q(\k(z),\Lambda)$ and $\omega\in\H^{q-1}(\k(z),\Lambda(-1))$.
\item
We now have a subgroup $\Gamma_{\D_*}^q\leq\H^q(\k(C))$ and a map 
\[\delta:\Gamma_{\D_*}^q(\k(C))\longrightarrow\Delta_{\D}^{q+1}(C)
\isom\bigoplus_{C_{(0)}}H^{q-1}(\k(z),-1)\]
induced by the map $\delta:\H^q(\k(C))\to\bigoplus_{C_{(0)}}\H_z^{q+1}(C)$ of 
\eqref{maplocseqs}.
This is not the residue map $\partial=\partial_1$ in the Kato complex
\eqref{katocomplex}.  For if $\alpha_C\in\Gamma_{\D_*}^q(\k(C))$
and $z\in\S$ is on $C_i\cap C_j$ then 
$\delta_z(\alpha_C)=\hat\delta_z\circ(\res_{\k(C)|K_{i,z}\times K_{j,z}})(\alpha_C)$,
and since $\alpha_C$ is in $\Gamma_{\D_*}^q$ we identify $\delta_z(\alpha_C)$
with $\omega:=\partial_z(\alpha_{C_i})$ by Lemma~\ref{subgroup}.
On the other hand, 
$\partial_z(\alpha_C)=\partial_z(\alpha_{C_i})+\partial_z(\alpha_{C_j})=2\omega$
by \eqref{partialz}.
Thus we have a commutative diagram
\[\xymatrix{
&\Gamma_{\D_*}^q(\k(C))\ar@{=}[d]\ar[r]^-\delta&\bigoplus_{C_{(0)}}\H^{q-1}(\k(z),-1)\ar[d]^\tau\\
&\Gamma_{\D_*}^q(\k(C))\ar[r]^-\partial&\bigoplus_{C_{(0)}}\H^{q-1}(\k(z),-1)}\]
where $\tau=\id$ if $z\not\in\S$ and
$\tau=2\cdot\id$ if $z\in\S$
\end{enumerate}

The next definition is motivated by Saito's \cite[Definition 2.1]{Sai85}.

\Definition\label{csdef}
Assume the setup of \eqref{setup}.

\begin{enumerate}[\rm i)]
\item
Define the group $\H^q_\cs(C)$ by the exact sequence
\[\xymatrix{
0\ar[r]&\H^q_{\cs}(C)\ar[r]&\H^q(C)\ar[r]&\H^q(\k(C))}\]
View $\H^q_\cs(C)$ as the subgroup of {\it completely split elements} of $\H^q(C)$.
\item
Let $\H^q_\cs(X)$ denote the preimage of $\H^q_\cs(C)$ in $\H^q(X)$
under the proper base change isomorphism $\H^q(X)\isim\H^q(C)$ of \eqref{maplocseqs}.
\item
Let $\H^q_\cs(V)$ and $\H^q_\cs(F)$ denote the images of 
$\H^q_\cs(X)$ in $\H^q(V)$ and $\H^q(F)$, respectively,
under the natural maps $\H^q(X)\to\H^q(V)$ and $\H^q(X)\to\H^q(F)$.
\item 
Let 
$\H^q_\ns(X),\;\;\H^q_\ns(V),\;\;\H^q_\ns(F)$
denote the corresponding ``nonsplit'' quotients $\H^q(-)/\H^q_\cs(-)$.
\end{enumerate}

\Lemma\label{cslemma}
We have $\H^q_\cs(V)=\ker(\H^q(V)\to\H^q(\k(C)))$,
and if $\alpha\in\H^q_\cs(V)$ then the image $\alpha_{F_D}$ of $\alpha$ in $\H^q(F_D)$
is zero for every prime divisor $D$ on $X$.
\rm

\begin{proof}
Let $K_X=\ker(\H^q(X)\to\H^q(V))$.
By \eqref{maplocseqs} we have a commutative diagram of exact rows
\[\xymatrix{
0\ar[r]&K_X\ar[r]\ar[d]&\H^q(X)\ar[d]^{\wr}\ar[r]&\H^q(V)\ar[d]^j\ar[r]^-\partial&\bigoplus_\D\H^{q-1}(D,-1)\ar[d]^h\\
0\ar[r]&\H^q_\cs(C)\ar[r]&\H^q(C)\ar[r]&\H^q(\k(C))\ar[r]^-\delta&\bigoplus_{C_{(0)}}\H^{q+1}_z(C)
}\]
The map $h$ is injective (by Lemma~\ref{subgroup}(a)), so
if $\alpha\in \ker(j)$ then $\partial(\alpha)=0$ by diagram chase, hence there
exists an element $\alpha_X\in\H^q(X)$ mapping to $\alpha$, and $\alpha_X\in\H^q_\cs(X)$ again by diagram chase,
since $\alpha\in\ker(j)$.
Therefore $\alpha\in\H^q_\cs(V)$ by Definition~\ref{csdef}.
Conversely if $\alpha\in\H^q_\cs(V)$ then $\alpha$ is in the image of some $\alpha_X\in\H^q_\cs(X)$,
and $\alpha$ is in $\ker(j)$ by diagram chase.  This proves the first statement.

Suppose $\alpha\in\H^q_\cs(V)$, $\alpha_X\in\H^q_\cs(X)$ is a preimage in $\H^q(X)$,
and $\alpha_C\in\H^q_\cs(C)$ is the image of $\alpha_X$ under the proper base change isomorphism.
To show $\alpha_{F_D}=0$ for any prime divisor $D$ on $X$ 
it suffices by \eqref{residues} to show that $\partial_D(\alpha)=0$ and $\alpha(D)=0$.

Since the map $\H^q(X)\to\H^q(F)$ factors through $\H^q(\O_{X,D})$ for any prime divisor $D$
on $X$, we have $\partial_D(\alpha)=0$ for all prime divisors $D$ on $X$ by \cite[3.6]{CT95}.
Since $\alpha_C\in\H^q(C)$ and $\alpha$ is the image of $\alpha_C$ via the inverse of
the proper base change map (composed with $\H^q(X)\to\H^q(V)$), $\alpha$ and $\alpha_C$ satisfy
the hypotheses of \cite[Theorem 4.9]{BT11},
and by \cite[Theorem 4.9(d)]{BT11}
$\alpha(D)=\inf_{\k(z)|\k(D)}(\alpha_C(z))$ for all horizontal prime divisors $D$ on $X$ lying over $z\in C$.  
Since the value $\alpha_C(z)$ is defined to be the image in $\H^q(\k(z))$
of the $z$-unramified element $(\alpha_C)_{\k(C)_z}$, and already $(\alpha_C)_{\k(C)}=0$, we conclude $\alpha(D)=0$
for all horizontal primes $D$.
It remains to compute the value of $\alpha$ at an irreducible component $C_i$ of $C$.
But $\alpha(C_i)=(\alpha_C)_{\k(C_i)}$ by \cite[Theorem 4.9(a)]{BT11}.
Since $\alpha_C\in\H^q_\cs(C)$ this shows $\alpha(C_i)=0$, as desired.
We conclude $\alpha_{F_D}=0$ for all prime divisors $D$ on $X$.

\end{proof}

\Remark
By Lemma~\ref{cslemma} the residue maps $\partial_D$ and the restriction maps $\res_{V|F_D}$
on $\H^q(V)$ with respect to
prime divisors $D$ of $X$ are well defined on $\H^q_\ns(V)$, and hence on $\H^q_\ns(F)$.
Thus just as for the ordinary groups
we will say a class $\alpha$ in $\H^q_\ns(V)$ or $\H^q_\ns(F)$ is {\it defined} on an open set $U$
if $\partial_D(\alpha)=0$ for all $D\subset U$.

\Theorem\label{newglue}
Assume \eqref{setup}.
Fix $\D$ on $X$ as in \eqref{distinguished},
let $\H^q_\ns(F)$ denote the nonsplit quotient of Definition~\ref{csdef}, and
let $\Gamma_{\D_*}^q(\k(C))\leq \H^q(\k(C))$ be the group of Definition~\ref{glueatz}.
Then for all $q\geq 0$ there is a homomorphism
\[\lambda=\lambda_\D:\Gamma_{\D_*}^q(\k(C))\longrightarrow\H_\ns^q(F)\]
that fits into a commutative diagram
$$
\xymatrix{
\H^q_\ns(F)\ar[r]^-\partial&\bigoplus_{\D}\H^{q-1}(\k(D),-1)\\
\Gamma_{\D_*}^q(\k(C))\ar[r]^-\delta\ar[u]^{\lambda}&\bigoplus_{C_{(0)}}
\H^{q-1}(\k(z),-1)\ar[u]_\inf
}
$$
where $\delta$ is the map of Remark~\ref{delta}(c) and $\partial$ is the residue map.
Let $\alpha=\lambda(\alpha_C)$, where $\alpha_C=(\alpha_{C_i})\in\Gamma_{\D_*}^q(\k(C))$ (so each $\alpha_{C_i}$ is in $\H^q(\k(C_i))$). 
Then:
\begin{enumerate}[\rm a)]
\item
$\alpha$ is defined on the generic points of $C$,
and $\alpha(C_i)=\alpha|_{\k(C_i)}=\alpha_{C_i}$.
\item
The ramification locus of $\alpha$ is contained in $\D$.
\item
If $D\in\D$ intersects $C$ at $z\in C_i$, $\pi_D\in\O_{X,z}$ is a local equation for $D$, 
and $\alpha_{C,z}=\alpha^\circ+(\bar\pi_D)\cdot\omega$ is the splitting over $K_{i,z}$
as in Lemma~\ref{keylemma} (where $\bar\pi_D\in K_{i,z}$ is the image of $\pi_D$), then 
over $F_D$ we have the Witt decomposition
\[
\res_{F|F_D}(\alpha)
=\inf_{\k(z)|\k(D)}(\alpha^\circ)+(\pi_D)\cdot\inf_{\k(z)|\k(D)}(\omega)
\]
\item
If $\alpha_C$ is unramified at a point $z\in C$,
then $\alpha$ is unramified at any horizontal prime divisor $D$ lying over $z$, and 
$\alpha(D)=\inf_{\k(z)|\k(D)}(\alpha_C(z))$.
\end{enumerate}
\rm

\begin{proof}
Set $V=X-\D$ as in \eqref{gluingframework}, and let
Let $\Gamma\df\im(\H^q_\ns(V)\to\H^q(\k(C)))=\im(\H^q(V)\to\H^q(\k(C)))$.
We claim the following commutative diagram has exact rows.
\begin{equation}
\xymatrix{
0\ar[r]&\H^q_\ns(X)\ar[r]\ar[d]^\wr&\H^q_\ns(V)\ar[r]^-\partial\ar[d]^\wr
&\bigoplus_{\D}\H^{q-1}(D,-1)\ar[r]\ar[d]^\wr&\H^{q+1}(X)\ar[d]^\wr\\
0\ar[r]&\H^q_\ns(C)\ar[r]\ar@{=}[d]&\Gamma\ar[r]^-{\delta}\ar[d]&
\bigoplus_{\D}\Delta_D^{q+1}
\ar[r]\ar[d]&\H^{q+1}(C)\ar@{=}[d]\\
0\ar[r]&\H^q_\ns(C)\ar[r]&\H^q(\k(C))
\ar[r]^-\delta&\bigoplus_{C_{(0)}}\H_z^{q+1}(C)
\ar[r]&\H^{q+1}(C)\\
}
\label{doubleladder1}
\end{equation}
where $\delta$ is defined in Remark~\ref{delta}(c).
Note the vertical map $\H^q_\ns(X)\to\H^q_\ns(C)$ is a well defined isomorphism
since the proper base change isomorphism $\H^q(X)\to\H^q(C)$ maps $\H^q_\cs(X)$ onto $\H^q_\cs(C)$,
the map $\H^q_\ns(V)\to\Gamma$ is an isomorphism by definition of $\Gamma$ and $\H^q_\ns(V)$,
and the map $\bigoplus_\D\H^{q-1}(D,-1)\to\bigoplus_\D\Delta_D^{q+1}=:\Delta_\D^{q+1}$ is an isomorphism by 
Lemma~\ref{subgroup}(a) and Definition~\ref{glueatz}.
Exactness of the top row at $\H^q_\ns(X)$ is immediate by \eqref{maplocseqs} since the map $\H^q(X)\to\H^q(V)$ takes 
$\H^q_\cs(X)$ onto $\H^q_\cs(V)$.
The residue map $\partial$ is well defined on $\H^q_\ns(V)$ by Lemma~\ref{cslemma}, 
and exactness at $\H^q_\ns(V)$ is immediate by the corresponding exactness at $\H^q(V)$ in
\eqref{maplocseqs}.
Exactness of the middle row follows from that of the top row since the two are isomorphic and the diagram commutes,
and exactness of the bottom row is immediate by \eqref{maplocseqs} and Definition~\ref{csdef}.

By Definition~\ref{glueatz}
we have $\Gamma_{\D_*}^q=\delta_\S^{-1}(\Delta_{\D_*}^{q+1})$, where $\Delta_{\D_*}^{q+1}=\bigoplus_{\D_*}\Delta_D^{q+1}$,
and $\delta_\S^{-1}(\Delta_{\D_*}^{q+1})=\delta^{-1}(\Delta_\D^{q+1})$ by Remark~\ref{delta}(a).
Therefore to show $\Gamma=\Gamma_{\D_*}^q$ it is enough to show
$\Gamma=\delta^{-1}(\Delta_\D^{q+1})$.
Suppose $\alpha_C\in\H^q(\k(C))$. 
An easy diagram chase using the exactness of the middle and bottom rows shows that 
if $\delta(\alpha_C)$ lies in the subgroup 
$=\Delta_\D^{q+1}\leq\bigoplus_{C_{(0)}}\H^{q+1}_z(C)$
then there exists a $\alpha_C'\in\Gamma$ such that $\delta(\alpha_C')=\delta(\alpha_C)$,
hence $\alpha_C'-\alpha_C$ is in the subgroup $\H^q_\ns(C)\leq\Gamma\leq\H^q(\k(C))$,
hence $\alpha_C\in\Gamma$.
Thus $\Gamma$ is precisely the set of elements in
$\H^q(\k(C))$ mapping to $\bigoplus_{\D}\Delta_D^{q+1}$, as desired.
Thus we have an isomorphism
\[\H^q_\ns(V)\overset{\sim}{\longrightarrow}\Gamma_{\D_*}^q\]
Composing the inverse maps $\Gamma_{\D_*}^q\to\H^q_\ns(V)$ and $\Delta_D^{q+1}\to\H^{q-1}(D,-1)$ 
with the maps $\H^q_\ns(V)\to\H^q_\ns(F)$ and $\H^{q-1}(D,-1)\to\H^{q-1}(\k(D),-1)$
yields the desired commutative square
$$
\xymatrix{
\H^q_\ns(F)\ar[r]^-\partial&\bigoplus_{\D}\H^{q-1}(\k(D),-1)\\
\Gamma_{\D_*}^q(\k(C))\ar[r]^-\delta\ar[u]^\lambda&\bigoplus_{\D}\Delta_D^{q+1}=:\Delta_\D^{q+1}(C)\ar[u]_\inf
}
$$

Suppose $\alpha_C=(\alpha_{C_i})\in\Gamma_{\D_*}^q(\k(C))$ and $\alpha=\lambda(\alpha_C)\in\H^q_\ns(F)$.
Then $\alpha$ comes from a unique element of $\H^q_\ns(V)$,
and is defined at the generic points of each $C_i$
since they are all contained in $V$.
Therefore as in \eqref{residues} $\alpha$ is in $\H^q_\ns(\O_{v_{C_i}})$,
and has a value $\alpha(C_i)=\alpha|_{\k(C_i)}$.
Since the map $\H^q(\k(C))\to\H^q_\ns(V)$ splits the restriction map $\res:\H^q_\ns(V)\to\H^q(\k(C))$, 
we have $\alpha|_{\k(C_i)}=\alpha_C|_{\k(C_i)}=\alpha_{C_i}$.
Therefore $\alpha(C_i)=\alpha|_{\k(C_i)}=\alpha_{C_i}$, proving (a).

Since $\alpha$ is in the image of $\H^q_\ns(V)\to\H^q_\ns(F)$, $\partial_{D'}(\alpha)=0$
for all $D'$ not in $\D$, therefore the ramification locus of $\alpha$ on $X$
is contained in $\D$, proving (b).
For $D\in\D$ running through $z$
the restriction $\res:\H^q(V)\to\H^q(F_D)$ factors through $\H^q_\ns(V)$
by Lemma~\ref{cslemma},
and also through $\H^q(\O_{X,z}^\h[1/\pi_D])$ by Lemma~\ref{keylemma},
at which point we have a decomposition $\alpha=\alpha^\circ+(\pi_D)\cdot\omega$
for some $\alpha^\circ\in\H^q(\k(z))$ and $\omega\in\H^{q-1}(\k(z),-1)$.
By Lemma~\ref{keylemma} we have the Witt decomposition
\[\alpha_{F_D}=\inf_{\k(z)|\k(D)}(\alpha^\circ)+(\pi_D)\cdot\inf_{\k(z)|\k(D)}(\omega)\]
Similarly the Witt decomposition with respect to $\bar\pi_D$
over $K_{i,z}$ is $\alpha^\circ+(\bar\pi_D)\cdot\omega$, proving (c).
The proof of (d)
proceeds exactly as in \cite[Theorem 4.9]{BT11}, and we omit the details.
This completes the proof.

\end{proof}

\Remark
Since we want to analyze $\H^q(F)$ and not $\H^q_\ns(F)$ we will often refer to lifts $\alpha\in\H^q(F)$
of some $\alpha_C\in\Gamma^q_{\D_*}(\k(C))$, by which we mean $\alpha$ is any preimage of $\lambda(\alpha_C)\in\H^q_\ns(F)$
under the projection $\H^q(F)\to\H^q_\ns(F)$.
By Lemma~\ref{cslemma}, $\alpha$ and $\lambda(\alpha_C)$ have identical images in $\H^q(F_D)$
for any prime divisor $D$ on $X$.

\Corollary\label{inreach}
In the situation of Theorem~\ref{newglue},
suppose $\alpha\in\H^q(F,\Lambda)$ has ramification divisor $D_\alpha$
that is horizontal and intersects each component of $C$ transversely.
Then there exists a choice of distinguished divisors $\D$
such that $\alpha$'s image in $\H^q_\ns(F,\Lambda)$ 
is in the image of $\lambda_{\D}:\Gamma_{\D_*}^q(\k(C),\Lambda)\to\H^q_\ns(F,\Lambda)$.
\rm

\begin{proof}
Since $D_\alpha$ intersects each component of $C$ transversely, 
we may assume that $D_\alpha\subset\D$, hence $\alpha$ is defined on $V=X-\D$.
By the proof of Theorem~\ref{newglue}, $\lambda$ is the inverse
of the natural isomorphism $\H_\ns^q(V)\to\Gamma_{\D_*}^q(\k(C))$, hence the corollary.

\end{proof}

\Remark\label{practicalglue}
To use Theorem~\ref{newglue} effectively we must be able to produce classes of $\H^q(\k(C))=\prod_i\H^q(\k(C_i))$ that are
in $\Gamma^q_{\D_*}(\k(C))$, and for this we must be able to construct elements of $\H^q(\k(C_i))$
with specified values in $\H^q(K_{i,z})$ for finitely many closed points $z\in C_i$.
This is possible in the following situations, which will come up in sections 3 and 4.

\begin{enumerate}[(a)]
\item
$\H^q(\k(C))=\H^1(\k(C),\Z/n)$ and each $\k(C_i)$ is the function field of a smooth curve.
Then we have Saltman's generalized Grunwald-Wang theorem \cite[Theorem 5.10]{Sa82},
which determines for any finite set of closed points $z_1,\dots,z_d$ on a given $C_i$,
and local data $\{\theta_{C_i,z_j}\in\H^1(K_{i,z_j},\Z/n)\}$, a global element $\theta_{C_i}\in\H^1(\k(C_i),\Z/n)$
with the given local data, provided we are not in the special case.
The {\it special case} for $\k(C_i)$, $n$, and the $K_{i,z_j}$
is the situation where $K_{i,z_j}(\mu_{2^{v_2(n)}})/K_{i,z_j}$ is noncyclic for some $j$.
Note that then necessarily $\car(k)=0$.
Thus there is no special case for $\k(C_i)$ if $k(\mu_{2^{v_2(n)}})/k$ is cyclic.
Even if there is a special case, if the $\theta_{C_i,z_j}$ are all {\it tractable} then a global $\theta_{C_i}$
still exists, by \cite[Theorem 5.8]{Sa82}.
\item
$\H^q(\k(C))=\H^1(\k(C),\mu_n)$ and each $\k(C_i)$ is the function field of a smooth curve.
Then we can use weak approximation as follows.
Fix a smooth $C_i$, closed points $z_1,\dots,z_d\in C_i$, and suppose given local data
$\{(g_j)\in\H^1(K_{i,z_j},\mu_n):j=1,\dots,d\}$.
By weak approximation there exists an element $f_i\in\k(C_i)$ such that $v_{z_j}(f_i-g_j)\geq n$
for each $z_j$.  
Then $f_i=g_j(1+\pi_j^{n-r}v)$ for some $r<n$, uniformizer $\pi_j$ for $z_j$, and $v\in\O_{C_i,z_j}$, 
and since the principal
units in $\O_{C_i,z_j}^\h$ are $n$-divisible (by Hensel's Lemma), the image of $(f_i)$ in $\H^1(K_{i,z_j},\mu_n)$
is $(g_j)$.
\item
$\H^q(\k(C))=\H^2(\k(C),\mu_n)$ and each $\k(C_i)$ is a global field of characteristic $p>0$.
Then we have the exact sequence
\[\xymatrix{
0\ar[r]&\H^2(\k(C_i),\mu_n)\ar[r]&\bigoplus_z\H^1(\k(z),\Z/n)\ar[r]&\H^1(k,\Z/n)\ar[r]&0}\]
Since $\H^1(\k(z),\Z/n)\isom\H^2(K_{i,z},\mu_n)$ this
sequence allows us to determine classes with specified local data.

\end{enumerate}

\section{Cyclic Classes in the Brauer Group}

\Paragraph
We assume the setup \eqref{setup} throughout this section.
We now apply Theorem~\ref{newglue} to study the structure of ${}_n\Br(F)=\H^2(F,\mu_n)$.
Our goal is to show that if $\alpha\in\H^2(F,\mu_n)$ is any class then there exists a model
$X$ for $F$ over which $\alpha$ decomposes into a sum of a Brauer element lifted from $\Gamma_{\D_*}^2(\k(C),\mu_n)$ and 
two $\Z/n$-cyclic cup products involving characters lifted from $\Gamma_{\D_*}^1(\k(C),\Z/n)$.
We then show that under a generalized Grunwald-Wang type hypothesis for algebras, $\H^2(F,\mu_n)$ is generated by $\Z/n$-cyclic
classes if the same is true for each $\H^2(\k(C_i),\mu_n)$, and that the $\Z/n$-length $\nL(F)$ is bounded in terms of $\nL(\k(C_i))$.
We make essential use not only of the above results, but also Saltman's
generalized Grunwald-Wang theorem from \cite[Theorem 5.10]{Sa82}, as in Remark~\ref{practicalglue}(a).
To apply this theorem we will need to hypothesize that 
{\it $k$ has no special case with respect to $n$}, meaning 
that $k(\mu_{2^{v_2(n)}})/k$ is cyclic.
The absence of an analogous theorem in higher degrees and for more general coefficients 
restricts our applications essentially to the Brauer group $\H^2(F,\mu_n)$.

The first problem is to compute the behavior of a given set of residues in $\bigoplus_{X_{(1)}}\H^1(\k(x),\Z/n)$
under blow-ups of $X$.
We use the notation of \eqref{katocomplex}.
Suppose $\psi\in\ker(\partial_1)\leq\bigoplus_{X_{(1)}}\H^1(\k(x),-1)$.
We will want to define a blow-up $\widetilde X\to X$ so that the 
transforms of the irreducible components of $\Supp(\psi)$
are regular on $\widetilde X$, and for this we need to see how $\psi$ pulls back to $\widetilde X$.
We write $\tilde\partial_q$ for the corresponding maps over $\widetilde X$,
then $\h_1(C_n^{-1}(\widetilde X))=\ker(\tilde\partial_1)/\im(\tilde\partial_2)$.

\Lemma\label{blowup}
Suppose $X$ is a connected regular noetherian surface.
Let $\widetilde X\to X$ be the blowup at a closed point $z$.  
Then there is a natural exact sequence
\[\xymatrix{0 \ar[r] &\H^1(\k(z),\Z/n) 
\ar[r] &\h_1(C_n^{-1}(\widetilde X)) \ar[r] &\h_1(C_n^{-1}(X)) \ar[r] &0}\]
\rm

\begin{proof}
Let $E$ be the special fiber of $\pi:\widetilde X\to X$ at $z$, 
and let $\xi_E$ be its generic point.
Consider the diagram of complexes
\begin{equation}
\label{blowupdiagram}
\xymatrix{
& & 0\ar[d] & 0\ar[d] & & \\
& 0\ar[r]\ar[d] & \H^1(\k(E),\Z/n)\ar[d]\ar[r]^-{\tilde\partial_1} &\ker(g)\ar[d]\ar[r] &0 \\
0\ar[r]&\H^2(F,\mu_n)\ar[r]^-{\tilde\partial_2}\ar@{=}[d]
&\bigoplus_{x\in \widetilde X_{(1)}}\H^1(\k(x),\Z/n)\ar[r]^-{\tilde\partial_1}\ar[d]^f
&\bigoplus_{y\in \widetilde X_{(0)}}\H^0(\k(y),\mu_n^*)\ar[d]^g\ar[r]&0\\
0\ar[r]&\H^2(F,\mu_n)\ar[r]^-{\partial_2}\ar[d]
&\bigoplus_{x\in X_{(1)}}\H^1(\k(x),\Z/n)\ar[r]^-{\partial_1}\ar[d]
&\bigoplus_{y\in X_{(0)}}\H^0(\k(y),\mu_n^*)\ar[r]\ar[d]&0\\
& 0 & 0 & 0 & & \\
}\end{equation}
where $f$ is the identity on 
$\H^1(\k(x),\Z/n)$ for $x\neq\xi_E$ and $f=0$ for $x=\xi_E$,
and $g$ is the identity on $\H^0(\k(y),\mu_n^*)$ for $y\in\widetilde X_{(0)}\backslash E$,
and on the complementary summand
\[g=\sum\cor:\bigoplus_{y\in E_{(0)}}
\H^0(\k(y),\mu_n^*)\longrightarrow\H^0(\k(z),\mu_n^*)\] 
With these maps the vertical sequences are all exact.

We claim each square of \eqref{blowupdiagram} is commutative.
This is obvious for the left half, so we will
show it only for the right two squares.
The map $\tilde\partial_1$
takes $\H^1(\k(E),\Z/n)$ into $\ker(g)$
by the exactness of the Faddeev sequence
\begin{equation}\label{Faddeevseq}
0\longrightarrow\H^1(\k(z),\Z/n) \longrightarrow\H^1(\k(E),\Z/n) \xrightarrow{\tilde\partial_1}\bigoplus_{y\in E_{(0)}} \H^0(\k(y),\mu_n^*)
\xrightarrow{\sum\cor}\H^0(\k(z),\mu_n^*) \longrightarrow 0
\end{equation}
hence the lower right square is commutative on elements of $\H^1(\k(E),\Z/n)$.

If $\tilde x\in\widetilde X_{(1)}$ and $\tilde x\neq\xi_E$
then $x=\pi(\tilde x)\in X_{(1)}$; let $\widetilde D$ and $D$ be the corresponding
divisors.
If $z\not\in D$ then $D$ and $\widetilde D$ have identical closed points
and $f$ is the identity on $\H^1(\k(\widetilde D),\Z/n)$, hence
$g\circ\tilde\partial_1=\partial_1$ on $\H^1(\k(\widetilde D),\Z/n)$.
If $z\in D$ then by \eqref{ramdivisor}
$\partial_{D,z}:\H^1(\k(D),\Z/n) \to \H^0(\k(z),\mu_n^*)$
is given by
\[\partial_{D, z}
= \sum_{v| z} \cor \circ \partial_{v}\]
where $v$ runs over all points of the normalisation of $D$ lying over
$z$, and $\partial_v: \H^1(\k(D),\Z/n) \to \H^0(\k(v),\mu_n^*)$ is the
residue map determined by $v$.  Since the
normalisation map factors through $\widetilde D \to D$, we have a
commutative diagram
\[\xymatrix{
\H^1(\k(\widetilde D),\Z/n)\ar@{=}[d]\ar[r]^-{\bigoplus \partial_{\tilde D,y}} 
& \bigoplus_{y| z} \H^0(\k(y),\mu_n^*)\ar[d]^{\sum\cor}\\
\H^1(\k(D),\Z/n)\ar[r]^-{\partial_{D, z}} & \H^0(\k(z),\mu_n^*)
}\]
where $y$ runs over all closed points of $\widetilde D$ lying over
$z$.  We conclude the lower right square is commutative, hence \eqref{blowupdiagram} commutes.

Write $A_\bullet$, $B_\bullet$ and $C_\bullet$ for the top, middle, 
and bottom rows of \eqref{blowupdiagram}, respectively.
Then we have an
exact sequence of complexes
\[\xymatrix{
0 \ar[r] &A_\bullet \ar[r] &B_\bullet \ar[r] &C_\bullet \ar[r] &0}\]
and associated long exact homology sequence
\[\xymatrix{
\cdots \ar[r] &\h_0(C_\bullet) \ar[r] &\h_1(A_\bullet) \ar[r] &\h_1(B_\bullet) \ar[r]
&\h_1(C_\bullet) \ar[r] &\h_2(A_\bullet) \ar[r] &\cdots
}\]
We claim $\h_0(C_\bullet)\to \h_1(A_\bullet)$ is the zero map.
For if $c\in\H^2(F,\mu_n)$ maps to zero in $\bigoplus\H^1(\k(D),\Z/n)$ then since $X$ is regular
$c$ is in the image of the map $\H^2(X,\mu_n)\to\H^2(F,\mu_n)$ by purity for (regular) surfaces
(see \cite[(4.2)]{Br10}).
It follows that the image of $c$ in $\H^1(\k(E),\Z/n)$ is zero, 
because $\H^2(X,\mu_n)\to\H^2(F,\mu_n)$ factors through $\H^2(\widetilde X,\mu_n)$,
and the map $\H^2(\widetilde X,\mu_n)\to\H^2(F,\mu_n)$ has image
\[\ker\left(\H^2(F,\mu_n)\to\H^1(\k(E),\Z/n)\oplus\bigoplus\H^1(\k(D),\Z/n)\right)\]
Therefore the image of $[c]\in \h_0(C_\bullet)$ in $\h_1(A_\bullet)$ is zero.
Also we have $\h_1(A_\bullet) =
\H^1(\k(z),\Z/n)$ and $\h_2(A_\bullet) = 0$ by the exactness of the
Faddeev sequence~(\ref{Faddeevseq}), and by definition $\h_1(B_\bullet)=\h_1(C_n^{-1}(\widetilde X))$,
and $\h^1(C_\bullet)=\h^1(C_n^{-1}(X))$.  Hence we obtain the desired exact
sequence
\[\xymatrix{
0 \ar[r] &\H^1(\k(z),\Z/n) \ar[r] &\h_1(C_n^{-1}(\widetilde X)) \ar[r] &\h_1(C_n^{-1}(X)) \ar[r] &0
}\]
\end{proof}

\Remark
It can be shown using Saltman's \cite[Theorem 5.2]{Sa08} that the sequence in Lemma~\ref{blowup}
is canonically split.  We will have no need for this result, so we omit the details.

\Corollary\label{preimage}
Assume \eqref{setup}.
If $\psi\in\ker(\partial_1)$,  then there exists a birational map
$\pi:\widetilde X\to X$ 
such that $\Supp(\tilde\psi)\cup\widetilde C$ has normal crossings,
where $\tilde\psi$ is any preimage of $\psi$ in $\ker(\tilde\partial_1)$, 
and $\widetilde C$ is the reduced subscheme of the
closed fiber of $\widetilde X$.
\rm

\begin{proof}
We can blow up $X$ until the preimage $\pi^*(\Supp(\psi)\cup C)$ has normal crossings
by Lipman's theorem on embedded resolution of singularities (see \cite[Theorem 2.26]{Liu}).
By Lemma~\ref{blowup} there is a preimage $\tilde\psi\in\ker(\tilde\partial_1)$ of $\psi$,
and the result follows since for any such $\tilde\psi$
the scheme $\Supp(\tilde\psi)\cup\widetilde C$ is a subset of $\pi^*(\Supp(\psi)\cup C)$.

\end{proof}

Recall from \eqref{ramdivisor} that if $\theta$ is an element of 
$\bigoplus_{\Div X}\H^1(\k(D),\Z/n)$ we write $\theta=(\theta_D)$, where $\theta_D$ is the $D$-th component,
and define the {\it support} $\Supp(\theta)$ of $\theta$ to be (the closed subset underlying) the effective divisor 
$\sum_{\theta_D\neq 0}D$, where the sum is over prime divisors on $X$.

\Lemma\label{2cyclics}
Assume \eqref{setup}, and that $k$ has no special case with respect to $n$.
Suppose $\psi\in\ker(\partial_1)$ over $X$.
Then there exists a birational map $\widetilde X\to X$ such that for any $\tilde\psi\in\ker(\tilde\partial_1)$
mapping to $\psi$, there exists an $\alpha\in\H^2(F,\mu_n)$ that is 
the sum of two $\Z/n$-cyclic classes,
such that $\tilde\psi-\tilde\partial_2(\alpha)$ is supported on a set of distinguished
divisors $\widetilde\D$ on $\widetilde X$.
\rm

\begin{proof}
We choose $\pi:\widetilde X\to X$ as in Corollary~\ref{preimage}, and let $\tilde\psi\in\ker(\tilde\partial_1)$
be any element mapping to $\psi$.
Then since $\Supp(\tilde\psi)\cup\widetilde C$ has normal crossings we may assume that 
the horizontal components of $\Supp(\tilde\psi)$ are in a set of distinguished divisors $\widetilde\D$ on $\widetilde X$,
and that they avoid the singular points of $\widetilde C$.
In general the $\tilde\psi_{C_i}$ will not glue across $\widetilde\S$ as in Definition~\ref{glueatz}, and for this
reason we blow up further until the dual graph of $\widetilde C$ is bipartite.
We now suppress the ``tilde'' notation, since everthing now happens over $\widetilde X$.
Since the dual graph of $C$ is bipartite the irreducible components of $C$ form two curves $C^+$ and $C^-$
such that $C^+\cap C^-=\S$, and the irreducible components of $C^+$ (resp. $C^-$) are disjoint.
Then the vertical part of $\Supp(\psi)$ consists of two parts, 
$C_\psi^+\subset C^+$ and $C_\psi^-\subset C^-$,
and the irreducible components of $C_\psi^+$ (resp. $C_\psi^-$) are disjoint.
By Corollary~\ref{picardcor} there exist elements
$(f^+),(f^-)\in\H^1(F,\mu_n)$ and divisors $H_i$ in $\D_\S$ such that 
\begin{align*}
\div(f^+)&=C_\psi^++H_1\pmod n\\
\div(f^-)&=C_\psi^-+H_2\pmod n
\end{align*}

Our strategy is to cup $(f^+)$ with the lift $\chi^+\in\H^1(F,\Z/n)$ 
of a character $\chi^+_C\in\Gamma^1_{\D_*}$ that matches $\psi$ on $C_\psi^+$,
subtract the residues of this cup product from $\psi$, and proceed similarly for $C_\psi^-$.
This will leave support on $\D$ only.

Define $\chi^+$ as follows:  For $C_i\subset C^+$ 
let $\chi^+_{C_i}=\psi_{C_i}\in\H^1(\k(C_i),\Z/n)$.
For all $C_j\subset C^-$ intersecting a $C_i\subset C^+$ there exists an element
$\chi^+_{C_j}\in\H^1(\k(C_j),\Z/n)$ such that
$\chi^+_{C_j,z}=\chi^+_{C_i,z}$ for all $z\in C_i\cap C_j$, 
by Saltman's generalized Grunwald-Wang theorem \cite[Theorem 5.10]{Sa82}.
(Note that there is no special case by hypothesis.)
If $C_j\subset C^-$ does not intersect $C^+$ we set $\chi_{C_j}^+=0$,
though this doesn't happen if $X$ is connected.
Thus we have defined characters $\chi_{C_k}^+$ on each irreducible component $C_k$ of $C$,
and by construction they glue across $\S$ along $\D$.
Hence we obtain an element $\chi_C^+\in\Gamma^1_{\D_*}$,
and by Theorem~\ref{newglue} there exists a $\chi^+\in\H^1(F,\Z/n)$ defined at the 
generic points of $C$ such
that $\chi^+|_{C_i}=\chi^+_{C_i}\in\H^1(\k(C_i),\Z/n)$ for all $C_i\subset C$, 
and $D_{\chi^+}\subset\D$.

Define $\chi^-$ similarly:  For $C_j\subset C_\psi^-$
let $\chi^-_{C_j}=\psi_{C_j}\in\H^1(\k(C_j),\Z/n)$
and for $C_j\subset C^-\backslash C_\psi^-$ let $\chi^-_{C_j}=0$;
for $C_i\subset C^+$ intersecting $C_j\subset C^-$
let $\chi^-_{C_i}$ be any character in $\H^1(\k(C_i),\Z/n)$ such that
$\chi^-_{C_i,z}=\chi^-_{C_j,z}$ for all $z\in C_j\cap C_i$;
and for $C_i\subset C^+$ not intersecting $C^-$ set $\chi_{C_i}^-=0$.
Thus we obtain an element $\chi_C^-\in\Gamma^1_{\D_*}$, hence an element
$\chi^-\in\H^1(F,\Z/n)$ by Theorem~\ref{newglue} with 
$\chi^-|_{C_i}=\chi^-_{C_i}\in\H^1(\k(C_i),\Z/n)$,
and $D_{\chi^-}$ supported in $\D$.

Let $\alpha^+=(f^+)\cdot\chi^+$, $\alpha^-=(f^-)\cdot\chi^-$, $\alpha=\alpha^++\alpha^-$,
and $\psi'=\psi-\partial_2(\alpha)$.
We compute $\psi_{C_i}'$ for the various $C_i\subset C$:
\begin{enumerate}
\item
If $C_i\subset C_\psi^+$, then $v_{C_i}(f^+)=1$, $v_{C_i}(f^-)=0$, 
$\chi^+_{C_i}=\psi_{C_i}$,
and $\partial_{C_i}(\chi^+)=\partial_{C_i}(\chi^-)=0$,
hence $\partial_{C_i}(\alpha^+)=\psi_{C_i}$ and $\partial_{C_i}(\alpha^-)=0$.
Therefore $\psi'_{C_i}=0$.
\item
Similarly if $C_j\subset C_\psi^-$
then $\partial_{C_j}(\alpha^+)=0$ and $\partial_{C_j}(\alpha^-)=\psi_{C_j}$,
so $\psi'_{C_j}=0$.
\item
If $C_i\subset C\backslash(C_\psi^+\cup C_\psi^-)$, then $\psi_{C_i}=0$ by definition,
and since neither $f^+,f^-,\chi^+$, nor $\chi^-$ ramify at $C_i$,
$\partial_{C_i}(\alpha)=0$, hence $\psi_{C_i}'=0$.
\end{enumerate}
This shows $\psi_{C_i}'=0$ for all vertical components $C_i\subset C$.
If $D$ is horizontal and irreducible and $D\not\in\D$, then again $\psi_D=0$, since we have
assumed the horizontal support of $\psi$ is in $\D_\S$, 
and then since $v_D(f^\pm)=0$ and $\partial_D(\chi^\pm)=0$, we have $\psi'_D=0$.
We conclude the support of $\psi'=\psi-\partial_2(\alpha)$ is in $\D$,
where $\alpha$ is the sum of two cyclic classes in $\H^2(F,\mu_n)$.

\end{proof}

\Corollary\label{plus2}
Assume \eqref{setup} and that $k$ has no special case with respect to $n$.
If $\alpha\in\H^2(F,\mu_n)$ then there exists a model $X/R$
such that 
\[\alpha=\beta+(f_1)\cdot\chi_1+(f_2)\cdot\chi_2\]
where $f_i\in F$, and the elements $\beta\in\H^2(F,\mu_n)$ and $\chi_1,\chi_2\in\H^1(F,\Z/n)$ are all lifted from $\k(C)$
as in Theorem~\ref{newglue}.
\rm

\begin{proof}
By Lemma~\ref{2cyclics} there exist cyclic classes $\gamma_1=(f_1)\cdot\chi_1$ and $\gamma_2=(f_2)\cdot\chi_2$
and a model $X$ with distinguished divisors $\D$
such that the divisor $D_\beta$ of the element $\beta=\alpha-\gamma_1-\gamma_2$ is supported on $\D$,
and the characters $\chi_1,\chi_2\in\H^1(F,\Z/n)$ are lifted from $\Gamma_{\D_*}^1$.
Since $D_\beta\subset X$ is horizontal and intersects each component of $C$ transversely,
$\beta$ is in the image of $\Gamma_{\D_*}^2$ by Corollary~\ref{inreach}.

\end{proof}

\Definition\label{cyclicdef}
Assume \eqref{setup}.  We say a $\Z/n$-cyclic class $\gamma_C\in\H^2(\k(C),\mu_n)$ 
is {\it cyclic along $\D_*$} ({\it via $\chi_C$}) if there exists
a $\chi_C\in\Gamma_{\D_*}^1(\k(C),\Z/n)\leq\H^1(\k(C),\Z/n)$ 
such that $\gamma_C=(f_C)\cdot\chi_C$ for some $(f_C)\in\H^1(\k(C),\mu_n)$.

\Remark
It is equivalent to say $\gamma_C$ is split by the cyclic extension associated to some 
$\chi_C\in\Gamma^1_{\D_*}(\k(C),\Z/n)$.
For if $\gamma_C\in\H^2(\k(C),\mu_n)$ is $\Z/n$-cyclic along $\D_*$ via $\chi_C\in\Gamma_{\D_*}^1(\k(C),\Z/n)$ 
then each $\gamma_{C_i}$ is split by the extension $\k(C_i)(\chi_{C_i})$, hence
$\gamma_C|_{\k(C)(\chi_C)}=0$.  Conversely if $\chi_C\in\Gamma_{\D_*}^1(\k(C),\Z/n)$
splits $\gamma_C$ then $\gamma_C$ is $\Z/n$-cyclic along $\D_*$ (via $\chi_C$).

We show elements of $\Gamma^2_{\D_*}(\k(C),\mu_n)$ that are $\Z/n$-cyclic along $\D_*$ lift to cyclic classes
in $\H^2(F)$.

\Lemma\label{cyclic2cyclic}
Assume \eqref{setup}.
If $\gamma_C\in\Gamma_{\D_*}^2\leq\H^2(\k(C),\mu_n)$ is cyclic along $\D_*$ via $\chi_C\in\Gamma_{\D_*}^1(\k(C),\Z/n)$,
then any lift $\gamma\in\H^2(F,\mu_n)$ of $\lambda_{\D}(\gamma_C)\in\H^2_\ns(F,\mu_n)$ is cyclic, 
via any lift $\chi\in\H^1(F,\Z/n)$ of $\lambda_{\D}(\chi_C)\in\H^1_\ns(F,\Z/n)$.
\rm

\begin{proof}  
The ramification divisor $D_\chi$ is independent of the choice of representative
$\chi$ of $\lambda_\D(\chi_C)$ by Lemma~\ref{cslemma},
and $D_\chi\subset\D$ by Theorem~\ref{newglue}(b).
Let $F(\chi)$ denote the field extension defined by $\chi$,
and let $\rho:Y\to X$ be the normalization of $X$ in $F(\chi)$.
Note that then $F(\chi)$ is the function field of $Y$.
Then since $X$ is regular and $n$ is prime-to-$p$, the finite separable field
extension $F(\chi)/F$ is tamely ramified along $D_\chi$,
and $Y\to X$ is a tamely ramified cover by \cite[Lemma 3.2]{BT11}.
Since $D_\chi\subset\D$ and $X$ satisfies the setup of \eqref{setup},
$Y$ is regular, $C_Y=C\times_X Y=Y_{0,\red}$, 
each irreducible component of $C_Y$ is regular, 
$\S_Y=\rho^{-1}(\S)$ is the set of singular points of $C_Y$, exactly two irreducible
components of $C_Y$ meet at each $z\in\S_Y$,
and the support of the irreducible components of $D_Y$ for $D\in\D$
generate a set of distinguished
divisors $\D_Y$ for $Y$, all by \cite[Lemma 3.1]{BT11}.
Thus $Y$ satisfies \eqref{setup}.

We show that the restriction $\H^2(\k(C),\mu_n)\to\H^2(\k(C_Y),\mu_n)$
takes $\Gamma_{\D_*}^2$ to $\Gamma_{\D_{Y_*}}^2$.
Let $y\in\S_Y$ be a point lying over $x\in\S$.
Assume $x\in C_i\cap C_j$, and let $B_i$ and $B_j$
be the irreducible components of $C_Y$ passing through $y$.
Let $E\in\D_{Y_*}$ and $D\in\D_*$ be the distinguished divisors passing through $y$ and $x$,
so that $E\subset D_Y$.
Let $\pi_D\in\O_{X,x}$ be a local equation for $D$, and let $\pi_j$ and $\pi_i$
be images of $\pi_D$ in $\k(C_i)$ and $\k(C_j)$, respectively.
Let $K_{i,y}=\Frac\O_{B_i,y}^\h$ and $K_{j,y}=\Frac\O_{B_j,y}^\h$,
the henselizations at $y$ on $B_i$ and $B_j$. 
Let $\pi_{B_i}$ and $\pi_{B_j}$ be local equations
for $B_i$ and $B_j$ in $\O_{Y,y}$.
Each is a regular prime meeting $E$ transversely at $y$.

Since $Y$ is regular $\O_{Y,y}$ is factorial, and we have a factorization
$\pi_D=u\pi_E^e$ for some number $e$, where $u\in\O_{Y,y}^*$ and $\pi_E$ is a local
equation for $E$.
Thus the image of $(\pi_j)$ in $\H^1(K_{i,y},\mu_n)$
Witt-decomposes as $(\bar u_j)+e\cdot(\pi_{j_y})$, where $\bar u_j$ is the image of $u$ in $\k(y)$
under the map $\O_{Y,y}\to\O_{Y,y}/(\pi_E,\pi_{B_i})=\k(y)$, and $\pi_{j_y}$
is the image of $\pi_E$ in $K_{i,y}$, which is a uniformizer since $E$ meets $B_i$
transversely.
Similarly the image of $\pi_i$ in $\H^1(K_{j,y},\mu_n)$ decomposes as $(\bar u_i)+e\cdot(\pi_{i_y})$,
where $\bar u_i$ is the image of $u$ in $\k(y)$
under the map $\O_{Y,y}\to\O_{Y,y}/(\pi_E,\pi_{B_j})=\k(y)$, and 
the image of $\pi_E$ in $K_{j,y}$ is a uniformizer $\pi_{i_y}$.
Since $E$ intersects both $B_i$ and $B_j$ transversely we have
$(\pi_E,\pi_{B_i})=(\pi_E,\pi_{B_j})$, hence $\bar u_j=\bar u_i\df\bar u$.

Now $\gamma_C\in\H^2(\k(C),\mu_n)$ has decompositions
$\gamma_{i,x}=\gamma_i^\circ+(\pi_j)\cdot\omega_i$ and
$\gamma_{j,x}=\gamma_j^\circ+(\pi_i)\cdot\omega_j$ in $\H^2(K_{i,x},\mu_n)$
and $\H^2(K_{j,x},\mu_n)$, respectively, as in Lemma~\ref{subgroup},
and since $\gamma_C\in\Gamma_{\D_*}^2$ we have $\gamma_i^\circ=\gamma_j^\circ\df\gamma^\circ$ and
$\omega_i=\omega_j\df\omega$.
Thus 
\begin{align*}
\res_{K_{i,x}|K_{i,y}}(\gamma_C)&=\gamma^\circ+(\pi_j)\cdot\omega=\gamma^\circ+(\bar u)\cdot\omega+(\pi_{j_y})\cdot e\omega\\
\res_{K_{j,x}|K_{j,y}}(\gamma_C)&=\gamma^\circ+(\pi_i)\cdot\omega=\gamma^\circ+(\bar u)\cdot\omega+(\pi_{i_y})\cdot e\omega\\
\end{align*}
and these are Witt decompositions in $\H^2(K_{i,y},\mu_n)$ and $\H^2(K_{j,y},\mu_n)$, respectively.
Therefore $\gamma_{C_Y}$ glues at $y$ along $E$ by Lemma~\ref{subgroup}.
Since this analysis applies to every point $y\in\S_Y$, we conclude $\gamma_{C_Y}\in\Gamma^2_{\D_{Y_*}}$.
Thus the restriction $\H^2(\k(C),\mu_n)\to\H^2(\k(C_Y),\mu_n)$
takes $\Gamma_{\D_*}^2$ to $\Gamma_{\D_{Y_*}}^2$, as desired.

Let $V=X-\D$, $V_Y=Y-\D_Y$,
$\lambda=\lambda_{\D}$, and $\lambda_Y=\lambda_{\D_Y}$.
Since the restriction sends $\Gamma_{\D_*}^2$ to $\Gamma_{\D_{Y_*}}^2$, by base change we have
a commutative diagram
\[\xymatrix{
\H^2(V,\mu_n)\ar[r]&\H^2(V_Y,\mu_n)\\
\Gamma_{\D_*}^2\ar[u]^{\lambda}_\wr\ar[r]&\Gamma_{\D_{Y_*}}^2\ar[u]_\wr^{\lambda_Y}
}\]
Thus if $(\gamma_C)_{\k(C_Y)}=0$, then $\gamma_{F(\chi)}=0$, hence $\gamma$ is cyclic via $\chi$.

\end{proof}

\Proposition\label{1cyclic}
Assume \eqref{setup} with $K=\Q_p$, so $F$ is the function field of a smooth $p$-adic curve.
Suppose $\gamma\in\H^2(F,\mu_n)$ has ramification divisor $D_\gamma$ on $X$ that
is horizontal and intersects each component of $C$ transversely.
Then $\gamma$ is cyclic.
\rm

\begin{proof}
By Corollary~\ref{inreach} there exists a model $X$, a choice of distinguished divisors
$\D$, and an element $\gamma_C\in\Gamma^2_{\D_*}(\k(C),\mu_n)$
such that $\gamma=\lambda_{\D}(\gamma_C)$.
It remains to show that $\gamma_C$ is cyclic by Lemma~\ref{cyclic2cyclic}.
Since each $\k(C_i)$ is a global field, the Hasse splitting principle applies in the Brauer
group, i.e., each factor $\gamma_{C_i}$
is split by any cyclic extension of $\k(C_i)$ that has local degree $n$ at each $\k(C_i)_z$
at which $\gamma_{C_i}$ is nontrivial.
We now easily construct a global $\chi_C\in\Gamma^1_{\D_*}(\k(C),\Z/n)$
such that $\k(C_i)(\chi_{C_i})$ splits $\gamma_{C_i}$ using the Grunwald-Wang Theorem
(note there is no special case here),
and this proves the Corollary.

\end{proof}

\Remark
By Corollary~\ref{plus2} any class in $\H^2(F,\mu_n)$ equals the sum of two cyclic classes
and an element lifted from $\H^2(\k(C),\mu_n)$ using Theorem~\ref{newglue}.  The proof
of Proposition~\ref{1cyclic} shows that if $F$ is the function field of a smooth
$p$-adic curve then the lifted class is cyclic, hence $\nL(F)\leq 3$.
We will improve on this number using a different construction, which appears as 
Theorem~\ref{splitram} below.

More generally when $K=K_v$ is any complete discretely valued field
we would like to prove that if each $\H^2(\k(C_i),\mu_n)$ is generated by $\Z/n$-cyclic classes
then the same holds for $F$, and we would like to bound $\nL(F)$ in terms of $\nL(\k(C_i))$.
For this we need the following hypothesis:

\Paragraph{\bf Hypothesis.}\label{hypothesis}
Let $\bar F$ be a field, $n\in\N$ prime-to-$\car(\bar F)$,
and $\alpha\in\H^2(\bar F,\mu_n)$.
Suppose 
$\nL(\alpha/\bar F)=N$,
$v_1,\dots,v_d$ are discrete valuations on $\bar F$, and
we have ``local data''
$\alpha_{v_i}
=\sum_{j=1}^N\gamma_{v_i}^{(j)}\in\H^2(\bar F_{v_i},\mu_n)$ 
where each $\gamma_{v_i}^{(j)}$ is $\Z/n$-cyclic.
Then there exist $\Z/n$-cyclic classes $\gamma^{(j)}\in\H^2(\bar F,\mu_n)$ such that 
$\alpha=\sum_{j=1}^N\gamma^{(j)}$
and
$\gamma_{\bar F_{v_i}}^{(j)}=\gamma_{v_i}^{(j)}$ for $j=1,\dots,N$.

This hypothesis holds, for example, if $\bar F$ is the function field of a smooth curve over a finite field,
by Remark~\ref{practicalglue}(c).

\Theorem\label{genbound}
Assume the setup of \eqref{setup}, $k$ has no special case with respect to $n$.
Suppose the function field $\bar F$ of any smooth curve over $k$ satisfies Hypothesis~\ref{hypothesis}
and that $\H^2(\bar F,\mu_n)$ is generated by $\Z/n$-cyclic classes.
Then $\H^2(F,\mu_n)$ is generated by $\Z/n$-cyclic classes, and if $\nL(\bar F)\leq N$ for all $\bar F$,
then $\nL(F)\leq 2(N+1)$.
\rm

\begin{proof}
By Corollary~\ref{plus2} and Proposition~\ref{cyclic2cyclic} it suffices to show that if 
$\alpha=\lambda(\alpha_C)\in\H^2(F,\mu_n)$ for
some $\alpha_C\in\Gamma^2_{\D_*}(\k(C),\mu_n)$ then we can write $\alpha_C=\sum_{j=1}^{P}\gamma_C^{(j)}$
for some $P\in\N$,
where $\gamma_C^{(j)}=(\gamma_{C_i}^{(j)})\in\Gamma_{\D_*}^2$ is $\Z/n$-cyclic
along $\D_*$,
and then that if $\nL(\k(C_i))\leq N$ for all $i$ 
then we can put $P\leq 2N$.
We will try to use Hypothesis~\ref{hypothesis} as little as possible.

As in the proof of Lemma~\ref{2cyclics} we may assume $C=C^+\sqcup C^-$ is bipartite.
Then since each $\H^2(\k(C_i))$ is generated by $\Z/n$-cyclic classes we have
$\nL(\alpha_{C_i}/\k(C_i))\leq M$ for all $C_i\subset C^+$, for some (finite) absolute bound $M$,
and we can write \[\alpha_{C_i}=\sum_{k=1}^M(f_{C_i}^{(k)})\cdot\theta_{C_i}^{(k)}\]
for $\theta_{C_i}^{(k)}\in\H^1(\k(C_i),\Z/n)$ and $(f_{C_i}^{(k)})\in\H^1(\k(C_i),\mu_n)$.
Since $k$ has no special case with respect to $n$,
by the generalized Grunwald-Wang theorem (see Remark~\ref{practicalglue}(a))
there exists an element $\theta_{C}^{(k)}\in\Gamma^1_{\D_*}(\k(C),\Z/n)$ 
such that $(\theta_{C}^{(k)})_{\k(C_i)}=\theta_{C_i}^{(k)}$.
The restrictions $\theta_{C_j}^{(k)}$ for $C_j\subset C^-$ glue to the $\theta_{C_i}^{(k)}$ at each $z\in\S$,
but otherwise they are unknown.

By the procedure outlined in Remark~\ref{practicalglue}(b)
there exists an element $(f_{C}^{(k)})\in\Gamma_{\D_*}^1(\k(C),\mu_n)$ whose $C_i$-components are $(f_{C_i}^{(k)})$.
Then \[\beta_C=\alpha_C-\sum_{k=1}^M(f_{C}^{(k)})\cdot\theta_{C}^{(k)}\] is supported on $C^-$ only, that is,
$\beta_{C^+}=0$.
Since $(f_C^{(k)})$ is in $\Gamma^1_{\D_*}(\k(C),\mu_n)$ and
$\theta_C^{(k)}$ is in $\Gamma^1_{\D_*}(\k(C),\Z/n)$ it is immediate
that each $(f_{C}^{(k)})\cdot\theta_{C}^{(k)}$ is in $\Gamma^2_{\D_*}(\k(C),\mu_n)$ and is $\Z/n$-cyclic along $\D_*$,
and hence that $\beta_C$ is in $\Gamma^2_{\D_*}(\k(C),\mu_n)$.

Since $\H^2(\k(C_j),\mu_n)$ is generated by $\Z/n$-cyclic classes for each $C_j\subset C^-$,
we can write each $\beta_{C_j}$ as a sum of $N$ $\Z/n$-cyclic classes for some (absolute) bound $N$.
Since $\beta_{C^+}=0$ and $\beta_C\in\Gamma^2_{\D_*}$ we have $\beta_{C_j,z}=0$ for each $C_j\subset C^-$ and $z\in\S\cap C_j$.
By Hypothesis~\ref{hypothesis} there exist $N$ $\Z/n$-cyclic classes $\varepsilon_{C_j}^{(l)}$ such that
$\beta_{C_j}=\sum_{l=1}^N\varepsilon_{C_j}^{(l)}$ and
$\varepsilon_{C_j,z}^{(l)}=0$ for each $z\in\S\cap C_j$.
Since each $\varepsilon_{C_j,z}^{(l)}=0$, each $\varepsilon_{C_j,z}^{(l)}$ glues with the zero class on $C_i\subset C^+$ across $z$
along $\D_*$, and we obtain a $\Z/n$-cyclic class $\varepsilon_C^{(l)}\in\Gamma^2_{\D_*}(\k(C),\mu_n)$
with restrictions $\varepsilon_{C_j}^{(l)}$ for $C_j\subset C^-$ and $\varepsilon_{C_i}^{(l)}=0$ for $C_i\subset C^+$.
By the generalized Grunwald-Wang theorem each $\varepsilon_C^{(l)}$ is $\Z/n$-cyclic along $\D_*$,
since any $\chi_C^{(l)}\in\Gamma^1_{\D_*}(\k(C),\Z/n)$ such that $\chi_{C_j}^{(l)}$ splits $\varepsilon_{C_j}^{(l)}$
for $C_j\subset C^-$
will split the whole class $\varepsilon_C^{(l)}$.
Now $\beta_C=\sum_{l=1}^N\varepsilon_C^{(l)}$ is a sum of $N$ classes in $\Gamma^2_{\D_*}(\k(C),\mu_n)$ 
that are $\Z/n$-cyclic along $\D_*$.
Altogether we have
\[\alpha_C=\sum_{k=1}^M(f_{C}^{(k)})\cdot\theta_{C}^{(k)}+\sum_{l=1}^N\varepsilon_C^{(l)}\]
This shows that $\alpha_C$ is expressible as a sum of classes of $\Gamma^2_{\D_*}(\k(C),\mu_n)$
that are $\Z/n$-cyclic along $\D_*$, as desired.  If additionally we have $\nL(\k(C_i))\leq N$ for all $\k(C_i)$,
i.e., $M=N$, then we obtain $\nL(\alpha/F)\leq 2N$, which completes the proof.

\end{proof}

\Corollary
Assume the setup of \eqref{setup} and $k$ has no special case with respect to $n$.
Suppose that the function field of any smooth curve over $k$
has $\Z/n$-length one.
Then $\H^2(F,\mu_n)$ is generated by $\Z/n$-cyclic classes, and $\nL(F)\leq 4$.
\rm

\begin{proof}
Hypothesis~\ref{hypothesis} holds trivially in this situation, so the result is immediate by
Theorem~\ref{genbound}.

\end{proof}

\section{Splitting Ramification in the Brauer Group}

This section is devoted to the following result.

\Theorem\label{splitram}
Assume \eqref{setup} and that $k$ has no special case with respect to $n$.
Suppose $\alpha\in\H^2(F,\mu_n)$.
Then there exists a $\Z/n$-cyclic class 
$\gamma\in\H^2(F,\mu_n)$ and a $\Z/n$-cyclic extension $L/F$
such that $(\alpha-\gamma)|_L$ is unramified, 
where $\gamma=(f)\cdot\theta$ and $L=F(\psi)$ are defined by an element $f\in F$ 
and characters $\theta,\psi\in\H^1(F,\Z/n)$ 
that are lifted from 
some $\theta_C,\psi_C\in\Gamma_{\D_*}^1$ over some regular model $X/R$.
\rm

\begin{proof}
As usual we may assume $D_\alpha\cup C$ has normal crossings on $X$, that 
$D_\alpha\subset C\cup\D_\S$, and that the dual graph of $C$ is bipartite.
Let $\chi=(\chi_D)=\partial_2(\alpha)\in\bigoplus_{D\subset D_\alpha}\H^1(\k(D),\Z/n)$ 
denote the residues of $\alpha$ on $X$.
Since the dual graph of $C$ is bipartite the irreducible components of $C$ form two curves $C^+$ and $C^-$
such that $C^+\cap C^-=\S$, and the irreducible components of $C^+$ (resp. $C^-$) are disjoint.

Recall from \eqref{distinguished} that if $z\in X$ is a closed point then $D_z$ denotes
the distinguished divisor passing through $z$.  
We will write $\D_\S^+$ (resp. $\D_\S^-$) for the subsets of $\D_\S$ that meet $C^+$ (resp. $C^-$),
so that $\D$ is the disjoint union \[\D=\D_*\sqcup\D_\S^+\sqcup\D_\S^-\]
We write $\D_\S^+(D)$ (resp. $\D_\S^-(D),\D_*(D)$) 
for the components of a divisor $D$ that are in $\D_\S^+$ (resp. $\D_\S^-,\D_*$).
By Corollary~\ref{picardcor} there exists an element
$f\in F^*$ such that 
\[\div(f)=C+\D(D_\alpha)+H_f\pmod n\]
where $H_f$ is a sum (possibly with multiplicities)
of distinguished divisors that avoid $\S$ and $\D(D_\alpha)$.
Note that $\D(D_\alpha)=\D_\S(D_\alpha)$.
If $z\in\S$ is in $C_i\cap C_j$, let $\pi_{C_i,z},\pi_{C_j,z}$ be local equations for $C_i$, $C_j$
in $\O_{X,z}$, respectively, such that
$f=-\pi_{C_i,z}\pi_{C_j,z}$ is the prime factorization of $f$ in $\O_{X,z}$.
We may assume then that $D_z=\Spec\O_{X,z}/(\pi_{C_i,z}+\pi_{C_j,z})$ is in $\D_*$ and $\pi_{D_z}=\pi_{C_i,z}+\pi_{C_j,z}$.
If $z\not\in\S$ and $z$ is a singular point of $\div(f)$ on $C_k\cap D_z$ then let $\pi_{C_k,z}$ and $\pi_{D_z}$
be local equations for $C_k$ and $D_z$ in $\O_{X,z}$ such that $f=-\pi_{C_k,z}\pi_{D_z}\in\O_{X,z}$.
If $z\in C_k\backslash(\D(D_\alpha)\cup H_f)$ let $\pi_{C_k,z}$ be the image of $f$ in $\O_{X,z}$.

The first step is to construct an element $(f)\cdot\theta\in\H^2(F,\mu_n)$ that matches the residues
of $\alpha$ on $C^+$ and the ``unramified part'' of the residues of $\alpha$ on $\D_\S^-$,
with zero residues at $H_f$.
The element $\theta\in\H^1(F,\Z/n)$ will be any preimage in $\H^1(F,\Z/n)$
of the nonsplit $\lambda_\D$-image of a carefully constructed
$\theta_C\in\H^1(\k(C),\Z/n)$.
Briefly, $\theta_C\in\Gamma_{\D_*}^1\leq\H^1(\k(C),\Z/n)$ comprises the components
$\theta_{C_i}=\chi_{C_i}$ for $C_i\subset C^+$, and components $\theta_{C_j}$ for $C_j\subset C^-$
that are lifted from local data
using the generalized Grunwald-Wang theorem so that they have certain residues at $\D_\S^-(D_\alpha)$
and they glue with the $\theta_{C_i}$ at $\S$ along $\D_*$.
Note again we have used that $k$ has no special case with respect to $n$.

Set $\theta_{C_i}=\chi_{C_i}$ for all $C_i\subset C^+$.
For $z\in C_k$ a closed point and $z\not\in\S$, let
$\chi_{D_z}=\chi_{D_z}^\circ+(\bar\pi_{C_k,z})\cdot\partial_z(\chi_{D_z})$
be the Witt decomposition of $\chi_{D_z}$ in 
$\H^1(\k(D_z),\Z/n)$ with respect to the image of $\pi_{C_k,z}$ in $\k(D_z)$.
If $z\in\S$ is on $C_i\cap C_j$ for $C_i\subset C^+$, let $\pi_j$ be the image of $\pi_{D_z}$ in $\k(C_i)_z$,
and let $\theta_{C_i,z}^\circ+(\pi_j)\cdot\partial_z(\theta_{C_i})$ 
denote the corresponding Witt decomposition of $\theta_{C_i,z}$ in $\H^1(\k(C_i)_z,\Z/n)$.
Now for $z$ a closed point on $C_j\subset C^-$, set
\begin{equation}\label{thetaCj}
\theta_{C_j,z}=\begin{cases}
\theta_{C_i,z}^\circ+(\pi_i)\cdot\partial_z(\theta_{C_i})&\text{ if $z\in\S$ is on $C_i\cap C_j$}\\
\chi_{D_z}^\circ &\text{ if $D_z\subset\D_\S^-(D_\alpha)$}\\
0 &\text{ if $z\in H_f$}
\end{cases}
\end{equation}
where $\pi_i$ is the image of $\pi_{D_z}$ in $\k(C_j)_z$.
For each $C_j\subset C^-$ let $\theta_{C_j}\in\H^1(\k(C_j),\Z/n)$ be an element with these local descriptions,
which exists by the generalized Grunwald-Wang theorem since $k$ has no special case with respect to $n$.
The $\theta_{C_j}$ and $\theta_{C_i}$ glue along $\S$ by Lemma~\ref{subgroup} 
to give an element $\theta_C\in\Gamma^1_{\D_*}$, and
we obtain an element $\lambda_{\D}(\theta_C)\in\H^1_\ns(F,\Z/n)$ by Theorem~\ref{newglue}.
Let $\theta\in\H^1(F,\Z/n)$ be any element mapping to $\lambda_{\D}(\theta_C)$.
By Lemma~\ref{cslemma} $\theta$ and $\lambda_\D(\theta_C)$ have the same residues and restriction to completions
at prime divisors $D$ of $X$. 
Note $D_\theta\subset\D$ by Theorem~\ref{newglue}(b), and using Theorem~\ref{newglue}(c)
we compute
$D_\theta\cap \D(D_\alpha)=D_\theta\cap\D_\S^+\subset\D_\S^+(D_\alpha)$, 
and $D_\theta\cap H_f=\varnothing$.
Note also that $\D_*(D_\theta)$ may be nonempty, as it consists of all components through $z\in\S$
in $C_i\cap C_j$ at which $\partial_z(\chi_{C_i})=-\partial_z(\chi_{C_j})$ is nonzero.

Let $\gamma=(f)\cdot\theta$.
Then $\Supp D_\gamma\subset \Supp(\div(f)+D_\theta)$, so to compute the residues of $\gamma$
we may restrict to $D\subset C+\D(D_\alpha)+ H_f+ D_\theta$.
We compute
\[\partial_D(\gamma)=\begin{cases}
\chi_{C_i} &\text{ if $D=C_i\subset C^+$}\\
\theta_{C_j} &\text{ if $D=C_j\subset C^-$}\\
\chi_{C_i,z}^\circ-(\bar\pi_{C_i,z})\cdot\partial_z(\chi_{C_i,z})
&\text{ if $D=D_z\subset\D_\S^+(D_\alpha)$}\\
\chi_{D_z}^\circ
&\text{ if $D=D_z\subset\D_\S^-(D_\alpha)$}\\
v_D(f)\cdot\chi_{C_i}(z)&\text{ if $D=D_z\subset \D_\S^+(H_f)$ and $z\in C_i$}\\
0 &\text{ if $D\subset\D_\S^-(H_f)$}\\
-(f)\cdot\partial_D(\theta)
&\text{ if $D\subset D_\theta\backslash\D(D_\alpha)$}\\
\end{cases}\]
based on the general cup product computation \eqref{cupresidue},
\[\partial_D(\gamma)=v_D(f)\cdot\theta_{F_D}-(f)\cdot\partial_D(\theta)+(-1)\cdot v_D(f)\cdot\partial_D(\theta)\]
All of these are obvious except for $D=D_z\subset\D_\S^+(D_\alpha)$ and $D=D_z\subset H_f$.  
If $D=D_z\subset\D_\S^+(D_\alpha)$ then
$\theta_{F_D}=\theta_{C_i,z}^\circ+(\pi_D)\cdot\partial_z(\theta_{C_i,z})$
by Theorem~\ref{newglue}(c), $v_D(f)=1$ by definition of $f$, 
and since the image of $f$ in $\O_{X,z}$ is $-\pi_{C_i,z}\pi_D$, 
$\partial_D(\theta)=\partial_z(\theta_{C_i,z})$,
and $\theta_{C_i,z}=\chi_{C_i,z}$,
\begin{align*}
\partial_D(\gamma)&=\theta_{F_D}-(f)\cdot\partial_D(\theta)+(-1)\cdot \partial_D(\theta)\\
&=\theta_{C_i,z}^\circ+(\pi_D)\cdot\partial_z(\theta_{C_i,z})-(-\pi_{C_i,z}\pi_D)\cdot\partial_z(\theta_{C_i,z})
+(-1)\cdot\partial_z(\theta_{C_i,z})\\
&=\chi_{C_i,z}^\circ-(\pi_{C_i,z})\cdot\partial_z(\chi_{C_i,z})\\
\end{align*}
Suppose $D=D_z\subset H_f$. 
Since $\theta$ is in the image of $\lambda$ we have 
$\partial_D(\theta)=\partial_z(\theta_{C_{k,z}})$ by Theorem~\ref{newglue}.
If $C_k\subset C^+$ then $\theta_{C_k,z}=\chi_{C_k,z}$ by construction of $\theta$,
and $\partial_z(\chi_{C_k,z})=0$ since $H_f$ avoids $\D(D_\alpha)$, by \eqref{sumtozero}.
If $C_k\subset C^-$ then $\theta_{C_k,z}=0$ by construction of $\theta$.
Thus $\partial_D(\theta)=0$ for $D\subset H_f$,
hence $\partial_D(\gamma)=v_D(f)\cdot\theta_{F_D}$.
By construction $\theta_{F_D}=0$ if $C_k\subset C^-$,
and $\theta_{F_D}=\inf_{\k(z)|\k(D)}(\chi_{C_k}(z))\in\H^1(D,\Z/n)$ if $C_k\subset C^+$,
by Theorem~\ref{newglue}(c,d).

Let $\beta=\alpha-\gamma$.
Then $D_\beta\subset D_\alpha\cup D_\gamma$, and we compute
\[\partial_D(\beta)=\begin{cases}
0 &\text{ if $D=C_i\subset C^+$}\\
\chi_{C_j}-\theta_{C_j}&\text{ if $D=C_j\subset C^-$}\\
\chi_{D_z}^\circ-\chi_{C_i,z}^\circ
&\text{ if $D=D_z\subset \D_\S^+(D_\alpha)$}\\
(\pi_{C_j,z})\cdot\partial_z(\chi_{D_z})
&\text{ if $D=D_z\subset\D_\S^-(D_\alpha)$}\\
-v_D(f)\cdot\chi_{C_i}(z)&\text{ if $D=D_z\subset\D_\S^+(H_f)$ and $z\in C_i$}\\
0 &\text{ if $D\in \D_\S^-(H_f)$}\\
(f)\cdot\partial_D(\theta)
&\text{ if $D\subset D_\theta\backslash\D(D_\alpha)$}\\
\end{cases}\]
All of these are obvious except for $D=D_z\subset\D_\S^+(D_\alpha)$,
where using the fact that $\partial_z(\chi_{C_i,z})+\partial_z(\chi_D)=0$ (by \eqref{sumtozero}) we compute
\begin{align*}
\partial_D(\beta)=\partial_D(\alpha)-\partial_D(\gamma)
&=\chi_D^\circ+(\pi_{C_i,z})\cdot\partial_z(\chi_D)
-\left(\chi_{C_i,z}^\circ-(\pi_{C_i,z})\cdot\partial_z(\chi_{C_i,z})\right)\\
&=\chi_D^\circ+(\pi_{C_i,z})\cdot\partial_z(\chi_{D})
-\left(\chi_{C_i,z}^\circ+(\pi_{C_i,z})\cdot\partial_z(\chi_D)\right)\\
&=\chi_D^\circ-\chi_{C_i,z}^\circ\\
\end{align*}

The next (last) step is to split the residues of $\beta$ with a cyclic extension
of $L/F$ constructed from a character $\psi_C\in\Gamma_{\D_*}^1\leq\H^1(\k(C),\Z/n)$.
Let $\psi_{C_j}=\chi_{C_j}-\theta_{C_j}$ for each $C_j\subset C^-$,
and for $C_i\subset C^+$ let $\psi_{C_i}$ be a generalized Grunwald-Wang lift such that
\[
\psi_{C_i,z}=\begin{cases}
\partial_{D_z}(\beta)
&\text{ if $D_z\subset\D^+_\S(D_\beta)$}\\
\psi_{C_j,z} &\text{ if $D_z\subset\D_*$}\\
\end{cases}
\]
Note $\partial_{D_z}(\beta)\in\H^1(\k(z),\Z/n)=\H^1(D_z,\Z/n)$ if $D_z\subset\D_\S^+(D_\alpha)$.
By ``$\psi_{C_i,z}=\psi_{C_j,z}$'' when $D_z\subset\D_*$ we mean with respect to the Witt decomposition determined 
by $D_z$, as in Lemma~\ref{subgroup}.
Then we get a $\psi_C\in\Gamma_{\D_*}^1$, hence a $\psi\in\H^1(F,\Z/n)$ mapping to $\lambda_{\D}(\psi_C)\in\H^1_\ns(F,\Z/n)$
by Theorem~\ref{newglue}.  Again $\psi$ and $\lambda_\D(\psi_C)$ have the same residues and completions at prime
divisors $D$ of $X$, by Lemma~\ref{cslemma}.

Let $L=F(\psi)$; we claim $\beta_L$ is unramified.
By construction the divisor $D_\psi$ is in $\D$.
Since $X$ is regular and $L/F$ is a finite separable field extension
that is tamely ramified (along $D_\psi$),
the normalization $Y$ of $X$ in $L$ is a (finite, flat, regular) tamely ramified cover of $X$, 
by \cite[Lemma 3.2]{BT11}.  
Thus if $E\subset Y$ is a prime divisor then $E$ lies over a prime divisor $D\subset X$, and
by \eqref{residues} we have a diagram
\[
\xymatrix{
&\H^2(L,\mu_n)\ar[r]^-{\partial_E}&\H^1(\k(E),\Z/n)\\
&\H^2(F,\mu_n)\ar[u]^\res\ar[r]_-{\partial_D}&\H^1(\k(D),\Z/n)\ar[u]_{e\cdot\res}\\
}
\]
where $e=e(E/D)=|\partial_D(\psi)|$ is the ramification index,
$e\cdot\psi_{F_D}$ is in the canonical subgroup $\H^1(\k(D),\Z/n)$ of $\H^1(F_D,\Z/n)$, 
and $\k(E)=\k(D)(e\cdot\psi_{F_D})$.
If $D\not\subset D_\beta$ then clearly $\partial_E(\beta_{L})=0$.
Suppose that $D\subset D_\beta$.
We've seen that $D_\beta\subset C^-+\D(D_\alpha)+\D_\S^+(H_f)+ D_\theta$, and
since $\D(D_\alpha)=\D_\S(D_\alpha)$ and $D_\theta\cap\D_\S^+\subset\D^+_\S(D_\alpha)$, 
\[D_\beta\subset C^-\cup\D_\S^+\cup\D_\S^-(D_\alpha)\cup\D_\S^-(D_\theta)\cup\D_*(D_\theta)\]
We treat these cases in order:

If $D=C_j\subset C^-$ then $\k(E)=\k(C_j)(\psi_{C_j})$,
and since $\partial_{C_j}(\beta)=\psi_{C_j}$, $\partial_E(\beta_{L})=0$.

If $D=D_z\subset\D^+_\S(D_\beta)$ then $\psi|_{F_{D}}=\partial_{D}(\beta)\in\H^1(D,\Z/n)$
by Theorem~\ref{newglue}(c), hence $\k(E)=\k(D)(\partial_{D}(\beta))$, hence $\partial_E(\beta_{L})=0$.

If $D=D_z\subset\D^+_\S\backslash\D^+_\S(D_\beta)$ then $\partial_E(\beta_L)=0$ since $\partial_D(\beta)=0$.

If $D=D_z\subset\D_\S^-(D_\alpha)$, then 
$e=|\partial_{D}(\psi)|=|\partial_z(\psi_{C_j})|=|\partial_z(\chi_{C_j})|$, since $\theta_{C_j,z}$
is unramified at $z$ when $D\subset\D_\S^-(D_\alpha)$ by \eqref{thetaCj}.
Since $\partial_z(\chi_{C_j})+\partial_z(\chi_{D})=0$ by \eqref{sumtozero},
we have $|\partial_z(\chi_{C_j})|=|\partial_z(\chi_{D})|$,
so multiplication by $e$ splits $\partial_D(\beta)=(\pi_{C_j,z})\cdot\partial_z(\chi_{D})$.
Therefore $\partial_E(\beta_{L})=0$.

If $D=D_z\subset \D_\S^-(D_\theta)\backslash\D_\S^-(D_\alpha)$
then $e=|\partial_{D}(\psi)|=|\partial_z(\psi_{C_j})|=|-\partial_z(\theta_{C_j})|=|-\partial_{D}(\theta)|$, since $\chi_{C_j}$ is unramified
at $z$.
This is enough to
split $\partial_D(\beta)=(f)\cdot\partial_{D}(\theta)$, hence $\partial_E(\beta_{L})=0$.

If $D=D_z\subset \D_*(D_\theta)$ 
then $\partial_z(\chi_{C_i})+\partial_z(\chi_{C_j})=0$ by \eqref{sumtozero}, and using 
$\partial_z(\theta_{C_j})=\partial_z(\theta_{C_i})=\partial_z(\chi_{C_i})$ (by definition) we compute
\[e=|\partial_{D}(\psi)|=|\partial_z(\psi_{C_j})|=|\partial_z(\chi_{C_j})-\partial_z(\theta_{C_j})
|=|\partial_z(\chi_{C_j})-\partial_z(\chi_{C_i})
|=|-2\partial_{D}(\theta)|\]
The image of $f$ in $\O_{X,z}$ is $-\pi_{C_i,z}\pi_{C_j,z}$ by definition, 
and since $\pi_{D}=\pi_{C_i,z}+\pi_{C_j,z}$
is a local equation for $D$ the image of $f$ in $\k(D)$ is 
$\bar\pi_{C_i,z}^2$.
Now we compute 
\[\partial_D(\beta)=(f)\cdot\partial_{D}(\theta)=2(\bar\pi_{C_i,z})\cdot\partial_D(\theta)\in\H^1(\k(D),\mu_n)\]
and since $|2(\bar\pi_{C_i,z})\cdot\partial_D(\theta)|$ divides 
$e=|-2\partial_D(\theta)|$, again $\partial_E(\beta_L)=0$.

Thus in all cases $\beta_L=(\alpha-\gamma)_L$ is unramified.
This completes the proof.

\end{proof}

\section{$\Z/n$-Length and Brauer Dimension}

Recall the {\it $\Z/n$-length} $\nL(F)$ of $F$ is the smallest number $c\in\N\cup\{\infty\}$ such that
any class in ${}_n\Br(F)$ can be expressed as a sum of $c$ $\Z/n$-cyclic classes.

\Theorem\label{lengthtwo}
Suppose $F$ is the function field of a $p$-adic curve,
and $n$ is prime-to-$p$.
Then $\nL(F)=2$.
\rm

\begin{proof}
Since $n$ is prime-to-$p$, ${}_n\Br(F)=\H^2(F,\mu_n)$.
Suppose $\alpha\in\H^2(F,\mu_n)$.
By Theorem~\ref{splitram} there exists a $\Z/n$-cyclic class $\gamma$
such that the residues of $\beta=\alpha-\gamma$ are split by a $\Z/n$-cyclic field
extension $F(\psi)$.
Since $C$ is a projective curve over a finite
field we have $\H^2_\nr(F(\psi),\mu_n)=0$ (see e.g. \cite[Theorem 4.5]{Br10}), hence $\beta_{F(\psi)}=0$.
Therefore there exists an element $\tau\in F$ such that $\alpha=\gamma+(\tau)\cdot\psi$.
This shows $\nL(F)\leq 2$.

To show $\nL(F)\geq 2$ it suffices to produce a class of period $n$ and index $n^2$,
since such a class cannot be $\Z/n$-cyclic.
We will construct one over $F_{C_1}$, 
where $C_1$ is an irreducible component of $C$, and then define it over $F$ using Lemma~\ref{subgroup}
and Theorem~\ref{newglue}.

By the Grunwald-Wang theorem (note there is no special case since $\car(k)>0$) there exists
an element $\chi_1\in\H^1(\k(C_1),\Z/n)$ of order $n$ such that $\chi_{1,z}=0$
for each $z\in\S\cap C_1$.
By Hasse's local-global principle
there exists an element $\alpha_1\in\H^2(\k(C_1),\mu_n)$ of order $n$ whose local invariant
has order $n$ at some closed point 
$x\not\in\S$ at which $\chi_{1,x}=0$, and whose local invariants
are zero at each $z\in\S\cap C_1$.
By Corollary~\ref{picardcor} there exists an element
$(\pi)\in\H^1(F,\mu_n)$ such that $\div(\pi)=C+H\pmod n$, 
where $H$ is a sum of distinguished (horizontal) divisors.
Then $\alpha_1+(\pi)\cdot\chi_1$ is an element of $\H^2(F_{C_1},\mu_n)$
of period $n$, and by Witt-Nakayama's index formula it has index
$|\chi_1|\cdot\ind((\alpha_1)_{\k(C_1)(\chi_1)})=n^2$.

Let $\chi_C\in\Gamma_{\D_*}^1$ and $\alpha_C\in\Gamma_{\D_*}^2$ be the elements whose restrictions to 
$\k(C_1)$ are $\chi_1$ and $\alpha_1$, respectively, and whose restrictions to 
the other $\k(C_i)$ are zero.
Let $\chi\in\H^1(F,\Z/n)$ be a lift of $\lambda_\D(\chi_C)$ and set $\alpha=\lambda_\D(\alpha_C)\in\H^2(F,\mu_n)$.
Then $\delta=\alpha+(\pi)\cdot\chi$ is in $\H^2(F,\mu_n)$, and
$\delta_{F_{C_1}}=\alpha_1+(\pi)\cdot\chi_1$ by Theorem~\ref{newglue}(a).
Since $\ind(\delta)\geq\ind(\delta_{F_{C_1}})$, we have $\ind(\delta)\geq n^2$,
and so $\delta$ cannot be $\Z/n$-cyclic, hence $\nL(F)\geq 2$, hence $\nL(F)=2$.

\end{proof}

\Remark
In the terminology of \cite{Sa07},
if $n$ is prime, the point $x\in C_1$ used to produce the class $\delta$
of period $n$ and index $n^2$ is a {\it hot point} for $\delta$, which for a class
of prime period $n$ is a ``tell'' for index equal to $n^2$ (and not $n$), by \cite[Theorem 5.2]{Sa07}.
The requirement for $x$ to be a hot point is that
$x\in C_1\cap D_x$ be a singular point of $D_\delta$,
both $\partial_{C_1}(\delta)$ and $\partial_{D_x}(\delta)$ be $x$-unramified,
and $\br{\partial_{C_1}(\delta)(x)}\neq\br{\partial_{D_x}(\delta)(x)}$.
We show $x$ is a hot point.  First we compute $\partial_{C_1}(\delta)=\chi_1$,
and since $\chi_{1,x}=0$ we have both that $\chi_1$ is $x$-unramified and $\chi_1(x)=0$.
On the other hand since $(\pi)$ and $\chi$ are unramified at $D_x$ we have
$\partial_{D_x}(\delta)=\partial_{D_x}(\alpha)$, and since $\alpha$ is a $\lambda$-lift
we compute $\partial_{D_x}(\alpha)=\inf_{\k(x)|\k(D_x)}(\partial_x(\alpha_1))$
by Theorem~\ref{newglue}(c).
Over the local field $K_{1,x}$ we have 
$\res_{\k(C_1)|K_{1,x}}(\alpha_1)=(\bar\pi_{D_x})\cdot\chi_x$ for
an (unramified) element $\chi_x\in\H^1(\k(x),\Z/n)\leq\H^1(K_{1,x},\Z/n)$ and $\bar\pi_{D_x}$ a uniformizer
for $K_{1,x}$, which is provided by the image of $\pi_{D_x}$ in $K_{1,x}$.
In the language of the proof, the element $\chi_x$ is the local invariant of $\alpha_1$ at $x$,
which is nonzero by definition.
Thus $\partial_{D_x}(\delta)=\chi_{D_x}:=\inf_{\k(x)|\k(D_x)}(\chi_x)$.
Clearly $\chi_{D_x}$ is $x$-unramified, and $\chi_{D_x}(x)\neq 0$.
Thus $\br{\chi_1(x)}\neq\br{\chi_{D_x}(x)}$, which proves $x$ is a hot point.

Recall the {\it $n$-Brauer dimension} $\nBrdim(F)$ of $F$ is the smallest number $d\in\N\cup\{\infty\}$
such that any class in ${}_n\Br(F)$ has index dividing $n^d$.
The following corollary was originally proved in \cite[Theorem 3.4]{Sa97}.

\Corollary\label{brdim2}
Suppose $F$ is the function field of a $p$-adic curve, and $n$ is 
prime-to-$p$.
Then $\nBrdim(F)=2$.
\rm

\begin{proof}
It is immediate that $\nBrdim(F)\leq 2$ by Theorem~\ref{lengthtwo}.
On the other hand the class of period $n$ and index at least $n^2$ constructed in the proof of
Theorem~\ref{lengthtwo} shows $\nBrdim(F)\geq 2$, hence we have equality.

\end{proof}

\Corollary\label{decomposable}
Suppose $F$ is the function field of a smooth $p$-adic curve, $n$ is
prime-to-$p$, and $D$ is an $F$-division algebra of period $n$ and index $n^2$.
Then $D$ is decomposable.
\rm

\begin{proof}
By Theorem~\ref{lengthtwo} there are two $\Z/n$-cyclic central simple algebras
$A_1$ and $A_2$ such that $\M_r(D)\isom A_1\otimes_F A_2$ for some number $r$. 
Since each $A_i$ has degree $n$ and $D$ has degree $n^2$, we conclude $r=1$,
hence $D$ is decomposable as $A_1\otimes_F A_2$. 

\end{proof}

\bibliographystyle{abbrv} 
\bibliography{hnx.bib}

\begin{thebibliography}{10}

\bibitem{Al36}
A.~Albert.
\newblock Simple algebras of degree {$p^e$} over a centrum of characteristic
  {$p$}.
\newblock {\em Trans. Amer. Math. Soc.}, 40(1):112--126, 1936.

\bibitem{Br10}
E.~Brussel.
\newblock On {S}altman's {$p$}-adic curves papers.
\newblock In {\em Quadratic forms, linear algebraic groups, and cohomology},
  volume~18 of {\em Dev. Math.}, pages 13--39. Springer, New York, 2010.

\bibitem{BMT}
E.~Brussel, K.~McKinnie, and E.~Tengan.
\newblock Indecomposable and noncrossed product division algebras over function
  fields of smooth $p$-adic curves.
\newblock {\em Adv. in Math.}, 226:4316--4337, 2011.

\bibitem{BT11b}
E.~Brussel and E.~Tengan.
\newblock Division algebras of prime period $\ell\neq p$ over function fields
  of $p$-adic curves.
\newblock {\em Israel J. Math. (to appear)}.

\bibitem{BT11}
E.~Brussel and E.~Tengan.
\newblock Formal constructions in the {B}rauer group of the function field of a
  $p$-adic curve.
\newblock {\em Trans. Amer. Math. Soc. (to appear)}.

\bibitem{CT95}
J.-L. Colliot-Th{\'e}l{\`e}ne.
\newblock Birational invariants, purity and the {G}ersten conjecture.
\newblock In {\em {$K$}-theory and algebraic geometry: connections with
  quadratic forms and division algebras ({S}anta {B}arbara, {CA}, 1992)},
  volume~58 of {\em Proc. Sympos. Pure Math.}, pages 1--64. Amer. Math. Soc.,
  Providence, RI, 1995.

\bibitem{Fuj02}
K.~Fujiwara.
\newblock A proof of the absolute purity conjecture (after {G}abber).
\newblock In {\em Algebraic geometry 2000, {A}zumino ({H}otaka)}, volume~36 of
  {\em Adv. Stud. Pure Math.}, pages 153--183. Math. Soc. Japan, Tokyo, 2002.

\bibitem{GMS}
S.~Garibaldi, A.~Merkurjev, and J.-P. Serre.
\newblock {\em Cohomological invariants in {G}alois cohomology}, volume~28 of
  {\em University Lecture Series}.
\newblock Amer.\ Math.\ Soc., 2003.

\bibitem{EGAI}
A.~Grothendieck.
\newblock \'{E}l\'ements de g\'eom\'etrie alg\'ebrique. {I}. {L}e langage des
  sch\'emas.
\newblock {\em Inst. Hautes \'Etudes Sci. Publ. Math.}, (4):228, 1960.

\bibitem{EGAIV:c}
A.~Grothendieck.
\newblock \'{E}l\'ements de g\'eom\'etrie alg\'ebrique. {IV}. \'{E}tude locale
  des sch\'emas et des morphismes de sch\'emas. {III}.
\newblock {\em Inst. Hautes \'Etudes Sci. Publ. Math.}, (28):255, 1966.

\bibitem{EGAIV:d}
A.~Grothendieck.
\newblock \'{E}l\'ements de g\'eom\'etrie alg\'ebrique. {IV}. \'{E}tude locale
  des sch\'emas et des morphismes de sch\'emas {IV}.
\newblock {\em Inst. Hautes \'Etudes Sci. Publ. Math.}, (32):361, 1967.

\bibitem{HHK09}
D.~Harbater, J.~Hartmann, and D.~Krashen.
\newblock Applications of patching to quadratic forms and central simple
  algebras.
\newblock {\em Invent. Math.}, 178:231--263, 2009.

\bibitem{Ka86}
K.~Kato.
\newblock A {H}asse principle for two-dimensional global fields.
\newblock {\em J. Reine Angew. Math.}, 366:142--183, 1986.
\newblock With an appendix by Jean-Louis Colliot-Th{\'e}l{\`e}ne.

\bibitem{Liu}
Q.~Liu.
\newblock {\em Algebraic Geometry and Arithmetic Curves}, volume~6 of {\em
  Oxford Graduate Texts in Mathematics}.
\newblock Oxford University Press, Oxford, 2002.
\newblock Translated from the French by Reinie Ern{\'e}, Oxford Science
  Publications.

\bibitem{Mat08}
E.~Matzri.
\newblock All dihedral algebras of degree 5 are cyclic.
\newblock {\em Proc. Amer. Math. Soc.}, 136:1925--1931, 2008.

\bibitem{Mer83}
A.~Merkurjev.
\newblock Brauer groups of fields.
\newblock {\em Comm. Alg.}, 11:2611--2624, 1983.

\bibitem{MS83}
A.~Merkurjev and A.~Suslin.
\newblock ${K}$-cohomology of {S}everi-{B}rauer varieties and the norm residue
  homomorphism.
\newblock {\em Math. USSR Izv.}, 21(2):307--340, 1983.

\bibitem{M}
J.~Milne.
\newblock {\em \'{E}tale {C}ohomology}.
\newblock Princeton University Press, 1980.

\bibitem{Sai85}
S.~Saito.
\newblock Class field theory for curves over local fields.
\newblock {\em J. Number Theory}, 21(1):44--80, 1985.

\bibitem{Sa97}
D.~Saltman.
\newblock Division algebras over $p$-adic curves.
\newblock {\em J. Ramanujan Math. Soc.}, 12:25--47, 1997.
\newblock see also the erratum \cite{Sa98} and survey \cite{Br10}.

\bibitem{Sa98}
D.~Saltman.
\newblock Correction to division algebras over $p$-adic curves.
\newblock {\em J. Ramanujan Math. Soc.}, 13:125--129, 1998.

\bibitem{Sa07}
D.~Saltman.
\newblock Cyclic algebras over $p$-adic curves.
\newblock {\em J. Algebra}, 314:817--843, 2007.

\bibitem{Sa82}
D.~J. Saltman.
\newblock Generic {G}alois extensions and problems in field theory.
\newblock {\em Adv. in Math.}, 43(3):250--283, 1982.

\bibitem{Sa08}
D.~J. Saltman.
\newblock Division algebras over surfaces.
\newblock {\em J. Algebra}, 320(4):1543--1585, 2008.

\bibitem{Sur10}
V.~Suresh.
\newblock Bounding the symbol length in the {G}alois cohomology of function
  fields of {$p$}-adic curves.
\newblock {\em Comment. Math. Helv.}, 85(2):337--346, 2010.

\bibitem{Tei37}
O.~Teichm\"{u}ller.
\newblock Zerfallende zyklische $p$-algebren.
\newblock {\em J. reine angew. Math.}, 176:157--160, 1937.

\bibitem{Tig84}
J.-P. Tignol.
\newblock On the length of decompositions of central simple algebras in tensor
  products of symbols.
\newblock In {\em Methods in ring theory ({A}ntwerp, 1983)}, volume 129 of {\em
  NATO Adv. Sci. Inst. Ser. C Math. Phys. Sci.}, pages 505--516. Reidel,
  Dordrecht, 1984.

\end{thebibliography}

\end{document}